# CORNER PERCOLATION ON $\mathbb{Z}^2$ AND THE SQUARE ROOT OF 17[1]

### By Gábor Pete

*Microsoft Research*

We consider a four-vertex model introduced by Bálint Tóth: a dependent bond percolation model on $\mathbb{Z}^2$ in which every edge is present with probability 1/2 and each vertex has exactly two incident edges, perpendicular to each other. We prove that all components are finite cycles almost surely, but the expected diameter of the cycle containing the origin is infinite. Moreover, we derive the following critical exponents: the tail probability $\mathbb{P}(\text{diameter of the cycle of the origin} > n) \approx n^{-\gamma}$ and the expectation $\mathbb{E}(\text{length of a typical cycle with diameter } n) \approx n^\delta$, with $\gamma = (5 - \sqrt{17})/4 = 0.219\ldots$ and $\delta = (\sqrt{17} + 1)/4 = 1.28\ldots$. The value of $\delta$ comes from a singular sixth order ODE, while the relation $\gamma + \delta = 3/2$ corresponds to the fact that the scaling limit of the natural height function in the model is the additive Brownian motion, whose level sets have Hausdorff dimension 3/2. We also include many open problems, for example, on the conformal invariance of certain linear entropy models.

**1. Introduction and results.** We consider the following *corner percolation* model introduced by Bálint Tóth. Intuitively, it is the maximal entropy probability measure on subsets of the edges of the lattice $\mathbb{Z}^2$ in which the set of edges incident to any given vertex is one of the four possible corners: $\ulcorner, \urcorner, \llcorner, \lrcorner$. However, this seemingly local constraint allows for only $2N$ bits of free choice to determine a configuration in an $N \times N$ square, as follows. Take two doubly infinite sequences, $\{\xi(n)\}_{n \in \mathbb{Z}}$ and $\{\eta(m)\}_{m \in \mathbb{Z}}$, of i.i.d. $+/-$ signs, each possibility having probability 1/2. If $\xi(n) = +1$, then we keep the "even" edges of $\mathbb{Z}^2$ along the vertical line $\{n\} \times \mathbb{Z}$, that is, the edges $\{(n, 2k), (n, 2k + 1)\}$ for all $k \in \mathbb{Z}$, and we delete the "odd" edges $\{(n, 2k - 1), (n, 2k)\}$. If $\xi(n) = -1$, we delete the even and keep the odd

Received May 2007; revised May 2007.

[1]Supported in part by NSF Grant DMS-02-4479 and OTKA (Hungarian National Foundation for Scientific Research) Grants T30074 and T049398.

*Key words and phrases.* Dependent percolation, dimer models, critical exponents, additive Brownian motion, simple random walk excursions, conformal invariance.

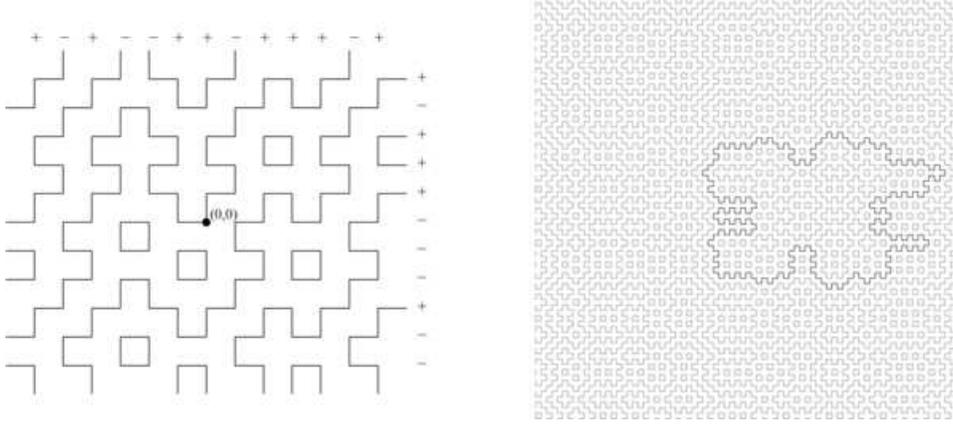

Fig. 1.  *Finite portions of some configurations.*

edges. Analogously, we delete every second edge along the horizontal line $\mathbb{Z} \times \{m\}$, according to the sign $\eta(m)$. Thus, we get a random 2-regular subgraph $G = G(\xi, \eta)$ of $\mathbb{Z}^2$; see Figure 1. Each connected component of $G$ is clearly a cycle or a bi-infinite path.

**Theorem 1.1.** *$G$ has no infinite components a.s. Moreover, each vertex is enclosed by infinitely many cycles. The expected diameter of the cycle containing the origin is infinite.*

Another way to look at this result is that a certain random walk on $\mathbb{Z}^2$ with long-term memory, which describes the component of the origin (and will be defined later on), is null recurrent. Given the fact that dependent percolation models and strongly self-interacting random walks are usually considered to be difficult, the proof of this first theorem is remarkably simple. The key step is to notice that looking at the components of $G$ as contour lines in a map, one can define a natural *height function* $H(n + \frac{1}{2}, m + \frac{1}{2})$ on the faces of the lattice $\mathbb{Z}^2$ (see Figure 3 below) which will turn out to equal $\lceil \frac{X_n + Y_m}{2} \rceil$, where $\{X_n\}_{n=-\infty}^{+\infty}$ and $\{Y_m\}_{m=-\infty}^{+\infty}$ are two independent bi-infinite simple random walks on $\mathbb{Z}$, certain functions of the sequences $\xi$ and $\eta$. Moreover, a detailed analysis of this connection makes it possible to understand the model quite thoroughly, for example, to determine the two most natural exponents.

**Theorem 1.2.** *The exponents for the tail probability*

$$\mathfrak{P}(n) := \mathbb{P}(\text{the diameter of the cycle of the origin } > n) \approx n^{-\gamma}$$

*and for the expectation*

$$\mathfrak{L}(n) := \mathbb{E}(\text{length of a typical cycle with diameter } n) \approx n^{\delta}$$

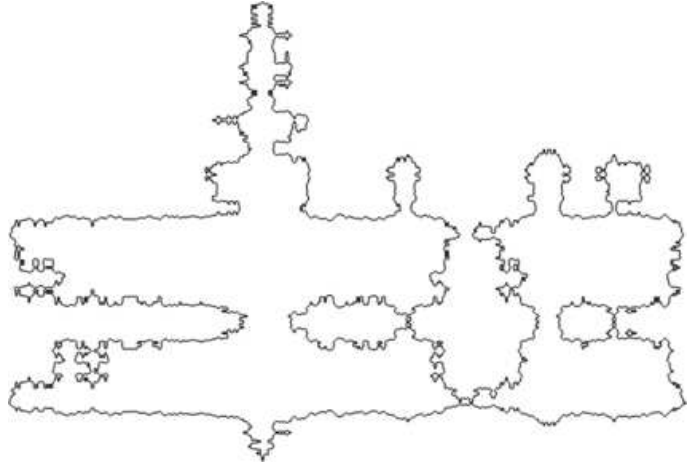

Fig. 2.   *A cycle of length 17,996.*

*exist, meaning that* $\gamma = \lim_{n \to \infty} \frac{-\log \mathfrak{P}(n)}{\log n}$ *and* $\delta = \lim_{n \to \infty} \frac{\log \mathfrak{L}(n)}{\log n}$. *The values are* $\gamma = (5 - \sqrt{17})/4 = 0.219\ldots$ *and* $\delta = (\sqrt{17} + 1)/4 = 1.28\ldots$.

We will define precisely what a "typical cycle" means in Section 4. The intuitive meaning of the exponent $\delta$ is the "dimension" of typical large cycles; see Figure 2. The relation $\gamma + \delta = 3/2$ corresponds to the fact that the height function $H(n + \frac{1}{2}, m + \frac{1}{2})$ has a natural scaling limit, the sum of two independent Brownian motions, $\mathcal{H}(t, s) = 1/2(W_t + W_s')$, and the level sets of this so-called *additive Brownian motion* have Hausdorff dimension 3/2 [36]. Finer geometric properties of these levels sets (which are often shared by those of the Brownian sheet) have been studied by Dalang, Mountford and Walsh in a series of papers [20, 21, 22, 23, 24, 25, 26, 44]; see [19] for a survey and, from a capacity point of view, for example, in [41]. Only after most of our work had been completed did we learn from Davar Khoshnevisan about the former papers, which contain the continuous analogs of some of our results (e.g., [22] finds a closed Jordan curve in the level set). One might guess that this curve should be the properly defined scaling limit of large cycles in our model, with Hausdorff dimension (which is presently unknown) $\delta = (\sqrt{17} + 1)/4 = 1.28\ldots$. However, passing to the limit seems to be a nontrivial task. See Section 7 for more details.

Another motivation for our model, which also provides some naive explanation as to why one would expect Theorem 1.1 to hold and the critical exponents $\gamma$ and $\delta$ to exist, is that corner percolation can be viewed as a linear entropy version of several classical conformally invariant models of statistical physics, as follows.

First, each edge of $\mathbb{Z}^2$ has probability $1/2$ to be in $G$, so we can view $G$ as *critical bond percolation* [33, 39] "conditioned" on the property that at each vertex, we see one of the four corners. However, this seemingly "local conditioning" introduces serious long-range dependence and ruins the global behavior of the model. Critical percolation is conformally invariant (proven at least on the triangular lattice by Smirnov [52]), while the additive Brownian motion appearing in our scaling limit is not. Furthermore, our tail exponent $\gamma$ for critical independent percolation is $5/48$ [43]; for other exponents, see [53]. Note that it is often the conformal invariance that is responsible for the rationality of certain critical exponents; in the physics literature, this is understood via connections to conformal field theory (CFT) [6], while, mathematically, it is understood via Schramm's stochastic Loewner evolutions (SLE) [48, 54]. See [3] for connections between SLE and CFT. Nevertheless, even without conformal invariance, the irrationality of such exponents is quite unusual.

Another conformally invariant process is the *double dimer model*, which is the union of two independent "uniform random perfect matchings" of $\mathbb{Z}^2$, while corner percolation is the union of a horizontal and a vertical perfect matching. In [37], Kenyon proved that the natural height function of the model has the *Gaussian free field* as its scaling limit, which is the conformally invariant two-time-dimensional version of Brownian motion; see [51]. Also, in the double dimer model, it is true that each vertex of $\mathbb{Z}^2$ is surrounded by infinitely many closed cycles almost surely [38, 50]. The height function fluctuations in an $n \times n$ subsquare of $\mathbb{Z}^2$ are of order $\log n$, as opposed to our $\sqrt{n}$, and it is conjectured that the exponent for the cycle length is $3/2$. A large range of similar random height function models is studied in [38] and [50], but those methods are not applicable in our strongly dependent model, where the Gibbsian "finite energy" condition completely fails.

Tóth defined corner percolation as a degenerate 4-vertex version of the famous *6-vertex model* [4, 5]. Many features of that model are not yet properly understood, but, being a generalization of the double dimer model, it is believed (sometimes proved—see, e.g., [31]) to behave in similar ways.

One can view the component of a fixed vertex in $G$ as the path of the following random walk on $\mathbb{Z}^2$ with long-term memory. At odd steps, we go up or down with probability $1/2$ each, at even steps, we go right or left with probability $1/2$ each, but whenever we visit one of the infinite horizontal or vertical lines we have already visited we must take the same direction as we did the first time on this line. If we consider the same walk without memory, just insisting on the alternating vertical/horizontal directions, then we get an "almost simple random walk," which is recurrent with infinite expected time of return, with the Brownian motion as its scaling limit. Thus, one might think, our path should also be closed, with infinite expected length. However, this argument is very weak. A little bit surprisingly, if we "interpolate" between these two recurrent walks (the almost simple random walk and the corner percolation path), by doing SRW (simple random walk) in the horizontal coordinate and the walk with memory vertically, then the resulting *random walk on a randomly oriented lattice* has been shown to be transient [14]. See Section 7 for more details on this model.

Since the distribution of the first version of this paper, corner percolation has inspired the definition of several linear entropy percolation models, including Benjamini's *trixor*, our *k-xor* and Angel's *odd-trixor* models on the triangular lattice; see Section 7 for definitions and some results. One shocking development is the following conjecture, strongly supported by computer simulations.

CONJECTURE 1.3. *The linear entropy $k$-xor ($k \geq 4$) and odd-trixor models on the triangular lattice have the same conformally invariant scaling limit as critical percolation.*

Corner percolation is also related to Winkler's nonoriented dependent percolation; see [57] and [2]. With a little additional work, which we will not present here, our methods imply the known result that there is no percolation in that model with three characters, while there is percolation in the case of four characters. However, there seems to be no natural meaning of our height function in that model; nor does the more puzzling oriented version [32] have a natural interpretation in our model.

Further aspects of the criticality of corner percolation, such as noise sensitivity/stability and the effect of biased coins, will also be discussed in Section 7, along with several more open problems.

We now turn to the discussion of the *height function*; see Figure 3. First, color the faces of $\mathbb{Z}^2$ black and white in a chessboard manner: let a face $(n + \frac{1}{2}, m + \frac{1}{2}) \in (\mathbb{Z} + 1/2) \times (\mathbb{Z} + 1/2)$ be black if $m + n$ is even and white if $m + n$ is odd. Fix the height of the black face touching the origin from the northeast to be $H(\frac{1}{2}, \frac{1}{2}) = 0$. Note that any component of $G(\xi, \eta)$, called a *contour line* hereafter, has all black faces along one side and all white faces along the other. Now, take a simple directed path through faces from $(\frac{1}{2}, \frac{1}{2})$ to some $(n + \frac{1}{2}, m + \frac{1}{2})$. Walking on this path, whenever we cross a contour from its black side toward its white side, we add 1 to the height and whenever we cross a contour from its white side toward its black one, we subtract 1. Clearly, the height $H(n + \frac{1}{2}, m + \frac{1}{2})$ determined in this way does not depend on the path taken. The *level of a contour* will be the height $H(\cdot, \cdot)$ of the black faces along the contour.

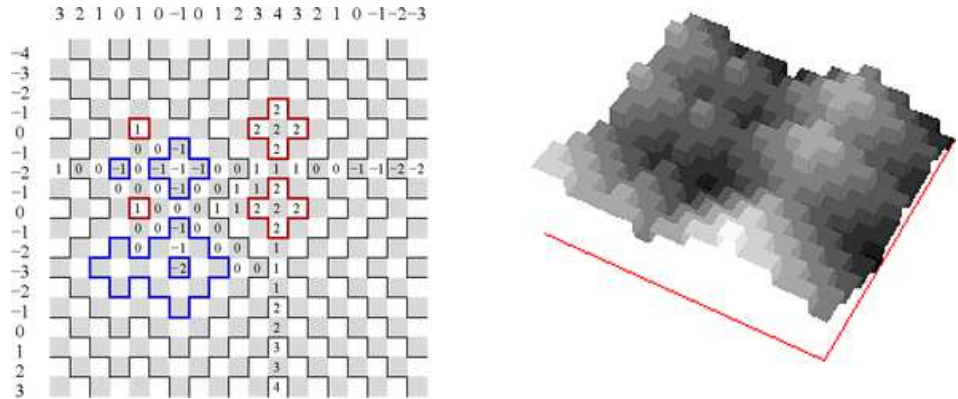

FIG. 3. *Definition of the height function.*

This height function can also be defined directly from the sign sequences $\{\xi(n)\}$ and $\{\eta(m)\}$. Let us first take $\xi^*(n) := (-1)^{n+1}\xi(n)$ and $\eta^*(m) := (-1)^{m+1}\eta(m)$. Then setting $X_0 = 0$, $X_n := \sum_{j=1}^{n} \xi^*(j)$ for $n > 0$ and $X_n := -\sum_{j=n+1}^{0} \xi^*(j)$ for $n < 0$, and similarly for $Y_m$, but with $\eta^*(j)$ instead of $\xi^*(j)$, defines two bi-infinite simple random walks, $\{X_n\}$ and $\{Y_m\}$, on $\mathbb{Z}$. It is now straightforward to check that

$$(1.1) \qquad H\left(n + \frac{1}{2}, m + \frac{1}{2}\right) = \left\lceil \frac{X_n + Y_m}{2} \right\rceil.$$

This single observation will be the key to our short proof of Theorem 1.1 in Section 2. Moreover, as we will see, most of the interesting phenomena in corner percolation can be formulated in a tractable way in terms of these simple random walks.

For example, having proven Theorem 1.1, we can divide all the contours in $G$ into two classes: a cycle is an *up-contour* if it has black faces along its exterior side; white faces along its interior side; it is a *down-contour* otherwise. In Figure 3, there are four up-contours and four down-contours; the directions of the other contours are impossible to determine from this finite piece of the configuration. It will turn out that each up-contour in $G$ can be identified with a so-called "compatible pair of up-excursions" in $\{X_n\}$ and $\{Y_m\}$, while down-contours are given by pairs of down-excursions. We will now define these compatible pairs.

An *up-excursion* of height $h \geq 1$ and length $2k \geq 2$ in $\{X_n\}$ is a subsequence $\{X_j\}_{j=a}^{a+2k}$, also denoted by $X[a, a + 2k]$, such that $X_a = X_{a+2k}$, while $X_j > X_a$ for all $a < j < a + 2k$ and $\max\{X_j - X_a : a \leq j \leq a + 2k\} = h$. Analogously, a *down-excursion* of height $h \geq 1$ and length $2k \geq 2$ in $\{X_n\}$ is a subsequence $X[a, a + 2k]$ such that $X_a = X_{a+2k}$, while $X_j < X_a$ for all $a < j < a + 2k$ and $\max\{X_a - X_j : a \leq j \leq a + 2k\} = h$. We say that two up-excursions (resp. two down-excursions) $X[a, a + 2k]$ and $Y[b, b + 2\ell]$ form a *compatible pair* if they have the same height $h$ and $X_a + Y_b + h$ is even (resp. odd).

The main combinatorial statement concerning the model is the following proposition, for which we will give a short proof in Section 3 using the height function representation. Our original proof involved a delicate induction and, in fact, the strong combinatorial structure of contours revealed by that proof (and well visible in Figure 2) was what led us to discover the representation by the simple random walks $\{X_n\}$, $\{Y_m\}$ and the height function (1.1) hidden in the original description of the model.

First, note that any closed cycle of $G$ has a smallest enclosing rectangle $[a+1,c] \times [b+1,d]$; the finite subsequences $X[a,c]$ and $Y[b,d]$ will be called the *marginals* of the rectangle.

**Proposition 1.4.** *The smallest enclosing rectangle of any up-contour (resp. down-contour) gives a pair of compatible up-excursions (resp. down-excursions) as marginals. Conversely, the rectangle given by any pair of compatible up-excursions (resp. down-excursions) has an up-contour (resp. down-contour) that spans this rectangle. The level of the contour is simply the value $X_a + Y_b + h$ given by the marginal excursions $X[a,c]$ and $Y[b,d]$. See Figures 3 and 5.*

With this tool available, we can now give a rough outline of a strategy to understand the exponents $\gamma$ and $\delta$. First, instead of insisting that the diameter be large or that it be exactly $n$, we will work with the height $h$ of the excursions that give the contour. Thus, we will consider the quantities $P(h) := \mathbb{P}(\text{the height of the marginal excursions of the cycle of the origin is } > h)$ and $L(h) := \mathbb{E}(\text{length of a cycle given by two compatible excursions of height } h)$. Since the length of an excursion of height $h$ divided by $h^2$ has a nontrivial limiting distribution as $h \to \infty$, with exponential tails, the statements

$$\gamma = \lim_{n \to \infty} \frac{-\log \mathfrak{P}(n)}{\log n}, \qquad \delta = \lim_{n \to \infty} \frac{\log \mathfrak{L}(n)}{\log n}$$

will turn out to be equivalent to

$$2\gamma = \lim_{h \to \infty} \frac{-\log P(h)}{\log h}, \qquad 2\delta = \lim_{h \to \infty} \frac{\log L(h)}{\log h}.$$

To compute $L(h)$, the main difficulty is that the actual shape of a contour is a rather strange "product" of the two marginal excursions, which will be described by the so-called "Two Cautious Hikers" algorithm in Section 3, and hence the length is difficult to find directly. It is much easier to give the total length of contours on that same level (say level 0) within the enclosing rectangle. This is roughly $ch^3$, corresponding to the fact that an $n \times n$ square has roughly $\sqrt{n}$ levels of the height function $H(\cdot, \cdot)$, with roughly $n^{3/2}$ faces on a typical level. However, a typical enclosing rectangle contains many more level-0 contours than just the one that touches its four sides (see, e.g., Figure 4). So, we need to subtract the total length of these additional contours, which all turn out to be contained entirely within our large enclosing rectangle (Lemma 3.1) and, hence, they are given by subexcursion pairs within our excursions. Such a subexcursion pair, with a certain height $m$, gives a contour of expected length $L(m)$. So, if we compute the expected number of compatible subexcursion pairs of a given height $m$ and giving a contour on level 0, then, multiplying this expected number by $L(m)$, we get their share in the total length of contours. Subtracting these products for all possible $m$'s, we can write a recursion for $L(h)$ involving all of the $L(m)$ with $m < h$. This plan will be carried out in Section 5, giving the recursion (5.6). This recursion can then be converted into a sixth order linear ordinary differential equation with nonconstant coefficients for the generating function of the sequence $L(h)$; see (5.18). This ODE has an isolated singularity of the first kind at 1 and the exact size of the singularity of the solution can be determined by the so-called *general Frobenius method* of linear ODE theory. The rate of growth of the sequence $L(h)$ can then be found using a standard Tauberian theorem.

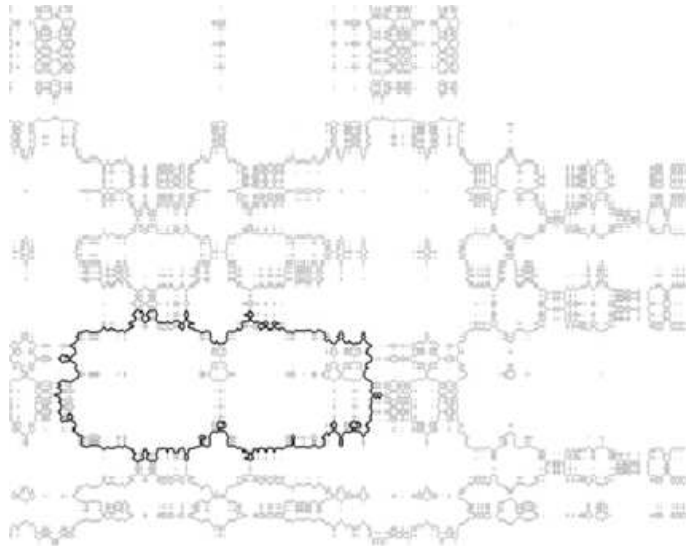

Fig. 4.   *The contours on a fixed level, with a distinguished large cycle.*

Writing a similar recursion for the sequence $P(h)$ seems harder, so we will find its decay rate by proving the relation $\gamma + \delta = 3/2$. This is done in Section 6 by considering the expected total number of edges on relatively large level-0 contour cycles in an $n \times n$ box.

**2. A short proof of Theorem 1.1.** We will prove a stronger statement: almost surely, for any fixed vertex and on any given level, there are infinitely many contour cycles surrounding that vertex.

Consider the square boundary

$$Q_N := \{(n+1/2, m+1/2) : \max\{|n|, |m|\} = N\}$$

in the dual lattice $(\mathbb{Z}+1/2) \times (\mathbb{Z}+1/2)$. Note that if, for some $0 < N < M$, the restrictions of the height function satisfy $H(Q_N) > 0$ and $H(Q_M) < 0$, then there must be a closed contour on level 0 surrounding the origin.

Now, consider the following event $\mathcal{A}_N = \mathcal{A}_N(\{X_n\})$ for a simple random walk $\{X_n\}_{n=-\infty}^{\infty}$ on $\mathbb{Z}$:

$$\mathcal{A}_N := \{X_n \geq -\sqrt{N} \text{ for } n \in \{-N, \ldots, N-1, N\}; X_{-N} > \sqrt{N}; X_N > \sqrt{N};$$

$$X_n \leq 2\sqrt{N} \text{ for } n \in \{-2N, \ldots, 2N-1, 2N\};$$

$$X_{-2N} < -2\sqrt{N}; X_{2N} < -2\sqrt{N}\}.$$

By Brownian scaling, there exists an absolute constant $\alpha > 0$ such that $\mathbb{P}(\mathcal{A}_N) > \alpha$ for all $N \geq 2$. Moreover, $\mathcal{A}_M$ becomes asymptotically independent of $\mathcal{A}_N$ as $M \to \infty$. In fact, there exists $0 < K < \infty$ such that $\mathbb{P}(\mathcal{A}_{KN}|\mathcal{A}_N)$ and $\mathbb{P}(\mathcal{A}_{KN}|\mathcal{A}_N^c)$ are both at least $\alpha/2$ for any $N > 0$.

The point of this construction is that if both the horizontal and vertical random walks, $\{X_n\}$ and $\{Y_n\}$, satisfy $\mathcal{A}_N$, then we have $H(Q_N) > 0$ and $H(Q_{2N}) < 0$. These two random walks are independent, thus $\mathcal{B}_N := \mathcal{A}_N(\{X_n\}) \cap \mathcal{A}_N(\{Y_n\})$ has $\mathbb{P}(\mathcal{B}_N) > \alpha^2$. Moreover, $\mathbb{P}(\mathcal{B}_{KN}|\mathcal{B}_N)$ and $\mathbb{P}(\mathcal{B}_{KN}|\mathcal{B}_N^c)$ are both at least $\alpha^2/2$ if $K$ is large enough, independently of $N$. Thus, in the sequence $\mathcal{B}_{K^i}, i = 1, 2, \ldots$, each event has a uniform positive probability to occur, independently of the occurrence or failure of the previous events. This means that infinitely many $\mathcal{B}_N$'s occur almost surely and we have infinitely many closed contours on level 0 surrounding the origin. Actually, this argument shows that $\mathbb{P}(\text{diameter} > n) \leq Cn^{-\beta}$, where $\beta > 0$ could be explicitly calculated.

To get infinitely many contours on an arbitrary level $\ell \in \mathbb{Z}$, we only have to consider translated versions of the events $\mathcal{A}_N$, which still have a uniform positive probability, so the same proof works.

To show that the expected diameter of the cycle of the origin is infinite, notice that whenever this cycle is not the smallest possible (of length 4), then it intersects the $x$- or the $y$-axis in a point different from the origin. So, with a fixed positive probability, this happens at the positive half of the $x$-axis, and at the intersection, there is a face with $H(N + \frac{1}{2}, \frac{1}{2}) = 0$. Since $Y_0 = 0$, this means that $X_N = 0$ for this $N$. If we take the smallest such $N > 0$, we get an excursion of length $N$ of the simple random walk $\{X_n\}$ on one hand and a chord of the contour cycle on the other. Since the probability that an SRW excursion has length larger than $n$ decays like $\pi^{-1/2}n^{-1/2}$ (see [29], Section 3.3), we have that $\mathbb{P}(\text{diameter} > n) \geq cn^{-1/2}$.

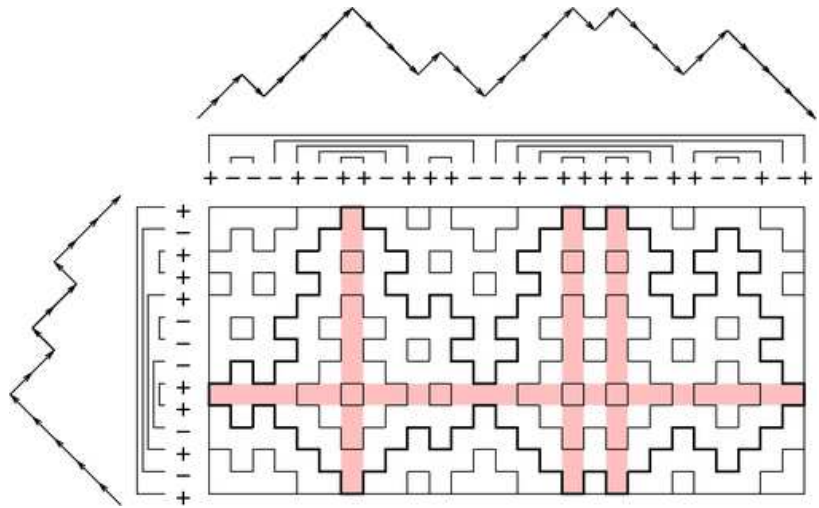

Fig. 5.   *A compatible pair determining a cycle with four passages.*

**3. The combinatorial description of contours.**   Before giving the proof of Proposition 1.4 using the height function representation, let us point out what the relevance of excursions in the $\{X_n\}$ and $\{Y_m\}$ sequences is to the configuration $G(\xi, \eta)$. Given the subconfigurations $\xi[a+1, c]$ and $\eta[b+1, d]$, replace the $k$th "+" sign with an opening parenthesis "(" if $k$ is odd and with a closing ")" if $k$ is even. Similarly, replace the $k$th "−" sign with an opening bracket "[" if $k$ is odd, and with a closing "]" if $k$ is even. Then, as it is immediate to see, we get a meaningful bracketing if and only if $X[a, c]$ and $Y[b, d]$ are both SRW excursions. As one can see from the proof below, if there is a contour cycle spanning the rectangle $\xi[a+1, c] \times \eta[b+1, d]$, then the corresponding matching between the opening and closing brackets and parentheses is "responsible" for the apparent symmetries of Figure 2. See also Figure 5.

PROOF OF PROPOSITION 1.4.   It will also follow from the proof that for each edge of a contour cycle that touches the enclosing rectangle, the opposite edge on the rectangle is also contained in the cycle and the row or column connecting these two edges is fully contained in the interior of the cycle. Such a row or column will be called a *passage*. We also claim that the passages occur exactly at the maxima of the marginal up-excursions and at the minima of the marginal down-excursions.

We will only consider up-contours and up-excursions. The case of down-contours and down-excursions can clearly be treated in the same way.

Suppose, first, that the walks $X[a, c]$ and $Y[b, d]$ form a compatible pair of up-excursions of height $h$, with $X_a = X_c = 0$ and $Y_b = Y_d = -h$. It is then clear that $X_n + Y_m > 0$ on the set

$$\mathcal{P} := \{(n,m) : X_n = h \text{ or } Y_n = 0; \ a < n < c \text{ and } b < m < d\},$$

but $X_n + Y_m \leq 0$ on the faces on the outer boundary of the rectangle $R :=$ $[a+1, c] \times [b+1, d]$. The set of faces given by $\mathcal{P}$ includes a cross that spans the rectangle $R$ in both directions, so the outer boundary of the component of $\{(n+1/2, m+1/2) : X_m + Y_m > 0\}$ containing this cross must be an up-contour spanning $R$.

Conversely, suppose that we have an up-contour, say on level 0, with enclosing rectangle $[a+1, c] \times [b+1, d]$. Then, for every $n \in [a, c]$, there is an $m' \in [b, d]$ with $X_n + Y_{m'} = 0$ and for every $n \in (a, c)$, there is an $m'' \in (b, d)$ with $X_n + Y_{m''} > 0$. Similarly, for every $m \in [b, d]$, there is an $n' \in [a, c]$ with $X_{n'} + Y_m = 0$ and for every $m \in (b, d)$, there is an $n'' \in (a, c)$ with $X_{n''} + Y_m > 0$. This implies that

$$\max X[a, c] - \min X[a, c] = \max Y[b, d] - \min Y[b, d]$$

and the minima can only be obtained at the endpoints. Moreover, there must be at least two places, $n_1, n_2 \in [b, d]$, with $(\max X[a, c]) + Y_{n_i} = 0$, thus $\min Y[b, d]$ must be actually obtained at both endpoints, so $Y[b, d]$ is an up-excursion. Similarly, $X[a, c]$ must also be an up-excursion. That they are compatible follows from the displayed equation.

Finally, it is clear that the set of faces corresponding to $\mathcal{P}$ is exactly the union of all the passages.  □

The following result concerning the global arrangement of contours on a fixed level will be used several times. See Figure 4.

LEMMA 3.1.    *If $C_i$, for $i = 1, 2$, are two contour cycles on the same level with enclosing rectangles $[a_i + 1, c_i] \times [b_i + 1, d_i]$, then one of the following three possibilities holds: $[a_1 + 1, c_1] \subseteq [a_2 + 1, c_2]$, $[a_2 + 1, c_2] \subseteq [a_1 + 1, c_1]$ or $[a_1 + 1, c_1] \cap [a_2 + 1, c_2] = \varnothing$. The analogous statement also holds for the vertical marginals. Furthermore, $C_2$ can never intersect the enclosing rectangle of $C_1$.*

PROOF.    If both cycles are up-contours, or both are down-contours, then the first statement is obvious from the possible ways that two up-excursions (or two down-excursions) can be arranged in one SRW trajectory. Note that, here, we did not even have to use the fact that the two contours are on the same level.

If $C_1$ is an up-contour and $C_2$ is a down-contour, then, by symmetry, it is enough to rule out the possibility that $a_1 < a_2 \leq c_1 < c_2$ and $b_1 < b_2 \leq d_1 < d_2$. Let us assume this possibility now. We may also assume that the level of the contours is 0 and so, for $m_i := \min X[a_i, c_i] = X(a_i) = X(c_i)$ and $M_i := \max X[a_i, c_i]$, and for $n_i := \min Y[b_1, d_1] = Y(b_1) = Y(d_1)$ and $N_i := \max Y[b_1, d_1] = N_1$, we have $m_1 + N_1 = n_1 + M_1 = 0$ and $m_2 + N_2 = n_2 + M_2 = 1$.

On the other hand, $a_1 < a_2 \le c_1$ implies $\max X[a_2, c_2] = X(a_2) = X(c_2) = M_2 \le M_1$, and $a_2 \le c_1 < c_2$ implies $\min X[a_2, c_2] = m_2 \le m_1$. Similarly, $\max Y[b_2, d_2] = Y(b_2) = Y(d_2) = N_2 \le N_1$ and $\min X[b_2, d_2] = n_2 \le n_1$. Combining all of these relations between $m_i, M_i, n_i$ and $N_i$, we find that $1 \le m_2 - m_1 \le 0$, a contradiction.

For the second statement, note that if $C_2$ intersected the enclosing rectangle $R_1$ of $C_1$, then, next to that intersection there would be a 0-valued face touching $R_1$ from the inside. But, by Proposition 1.4, this face would not only have a side which is an edge of $C_2$, but also another side which is an edge of $C_1$. From the existence of this face, one can easily see that the only way to have two distinct cycles, $C_1$ and $C_2$, would be if $C_1$ were a cycle of length 4; but, then, $R_1 = C_1$ and $C_2$ cannot intersect $R_1$ at all.  □

Though we will not use the following result later, let us describe how the exact shape of the contour cycle is determined by the compatible pair of excursions of height $h$.

Let the two up-excursions be $X[0, 2k]$ and $Y[0, 2\ell]$, with minima $X(0) = X(2k) = 0$, $Y(0) = Y(2\ell) = -h$ and maxima $X(m_i) = h$, $Y(n_j) = 0$, where the maximum places are $m_1 < \cdots < m_s$ and $n_1 < \cdots < n_t$, respectively. Consider the parts $X[0, m_1]$ and $Y[0, n_1]$, and let $Y'(j) := -Y(n_1 - j)$ for $j = 0, \ldots, n_1$. Imagine that two hikers, starting from a lake at the foot of two mountains that are described by the curves $X[0, m_1]$ and $Y'[0, n_1]$, want to reach the two peaks, one hiker for each peak, in such a way that their height levels always agree during hiking. It is known from a popular combinatorics exercise (at least in Hungary) that this is always possible: Consider the graph that has vertex set $V := \{(i, j) : X(i) = Y'(j), \ 0 \le i \le m_1, \ 0 \le j \le n_1\}$ and edge set $E$ connecting the vertices of $V$ that are accessible from each other in one step by the two hikers. Then every vertex has degree 2 or 4, except for $(0, 0)$ and $(m_1, n_1)$, which have degree 1. Hence, it is clear that there must be a path in $(V, E)$ connecting $(0, 0)$ and $(m_1, n_1)$, and this gives a path that the hikers can follow.

Of course, there can be several different paths between $(0, 0)$ and $(m_1, n_1)$. The following algorithm, which we call the algorithm of the *Two Cautious Hikers*, distinguishes one particular path, which will correspond to the part of the contour cycle that is determined by the excursion parts $X[0, m_1]$ and $Y[0, n_1]$.

Let us suppose that Xavier, the hiker climbing the mountain given by $X[0, m_1]$, is "cautious" when going up, while Yvonne, the hiker climbing the $Y'[0, n_1]$-mountain, is "cautious" when going down. This means that during their paths, it never happens that both of them backtrack, but whenever they arrive at a vertex of $V$ from which there are four edges going out, and at which they have never been before, they choose the edge on which Xavier backtracks if they have to continue upward, and on which Yvonne backtracks if they continue downward. When they arrive at a degree 4 vertex where they have already been, they choose the one edge which they have not previously used; see Figure 6. Starting from $(0, 0)$, this rule always gives a unique way to continue the path, until they reach $(m_1, n_1)$.

This path $(x_i, y_i)_{i=0}^{T}$ through the vertices of $V$, with $(x_0, y_0) = (0, 0)$ and $(x_T, y_T) = (m_1, n_1)$, gives a sequence of zeros $X(x_i) + Y(y_i) = 0$, while one of the two sums, $X(x_i) + Y(y_{i+1})$ and $X(x_{i+1}) + Y(y_i)$, is 1 and the other is $-1$, depending on whether the hikers went up or down from step $i$ to step $i + 1$. Mark an edge of $\mathbb{Z}^2$ if it is between such a 0- and a 1-valued face along the path. It is immediate to check that the marked edges form the part of the contour cycle inside the rectangle given by $X[0, m_1]$ and $Y[0, n_1]$, and the length of this part of the contour is just $2T$. The other three parts of the contour inside the rectangles given by the horizontal marginals $X[0, m_1]$, $X[m_s, 2k]$ and the vertical marginals $Y[0, n_1]$, $Y[n_t, 2\ell]$ can be described in the same way. The parts of the contour given by marginals $X[m_i, m_{i+1}]$ and $Y[n_j, n_{j+1}]$ can be handled by breaking these subexcursions into smaller mountains. Summarizing, we have proven the following result.

PROPOSITION 3.2. *The contour cycle given by a compatible pair of excursions can be naturally broken into several pieces, each of which corresponds to a path through the zero set $X_n + Y_m = 0$ that is given by the Two Cautious Hikers algorithm.*

**4. What is a typical contour cycle?** The natural definition for a *typical contour of diameter $n$* is the following. Consider corner percolation in a large $N \times N$ box, with $N \gg n$, and pick one of the diameter-$n$ cycles contained in the box, uniformly at random. Then, as $N \to \infty$, the distribution of this random cycle should converge to the distribution we seek to define. However, it would be very hard to work directly with this definition because, as we see from Proposition 1.4, the natural parameter for a cycle is not its diameter, but the height of the compatible excursions giving rise to it. Thus, we will consider "typical cycles of height $h$," and although our results will not imply that the large-$N$ limit distribution for the cycles of diameter $n$ exists, it will be clear that for large enough $N$, the expectation of the length of the randomly chosen cycle of diameter $n$ is close to $n^\delta$, in the sense given in Theorem 1.2.

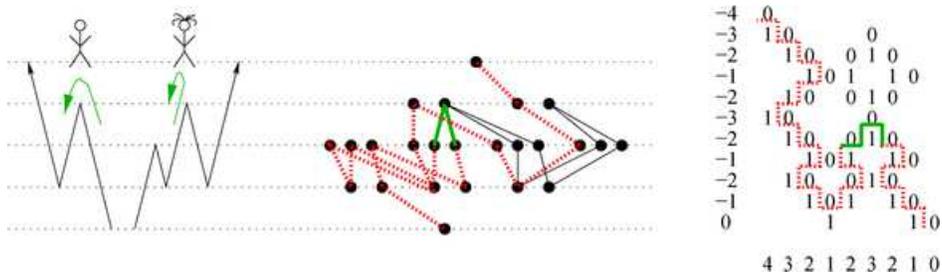

FIG. 6. *Part of a contour described by the Two Cautious Hikers.*

Similarly to the above, a *typical cycle of height $h$* is the large $N$ limit of a uniformly chosen cycle of height $h$ from the $N \times N$ box. But, this time, we have a direct way to describe this limit. By the symmetry between up- and down-contours and between different levels, and by simple properties of SRW excursions that we will discuss below, choosing one of these cycles uniformly at random from the box gives the same distribution on lattice cycles as choosing uniformly at random one of the height-$h$ up-excursions from the marginal random walk $X$ along the box and then, conditioned to have at least one compatible up-excursion in the $Y$ marginal (the event of which has a probability tending to 1 as $N \to \infty$), choosing one of them uniformly at random and considering the up-contour determined by this compatible pair. And the large-$N$ limit of this distribution is simply the up-contour determined by two independent up-excursions of height $h$. The simplest construction of an *excursion conditioned to have height $h$*, denoted by $\mathcal{E}_h$, is the following.

Run a simple random walk from 1, conditioned to reach $h$ before 0, until it reaches $h$. Then run another simple random walk, now started from $h$, conditioned to hit 0 before hitting $h + 1$, until it reaches 0. The two legs of this walk will be referred to as *there* and *back*. To compute anything about these conditioned walks, we will use the following basic lemma. This construction is a version of Doob's $h$-transform [28], with $h(i) := \mathbb{P}_i(\mathcal{A})$ below; we include the proof for completeness.

LEMMA 4.1.   *Let $\{X_n\}_{n=0}^\infty$ be any time-homogeneous Markov chain on the state space $\mathbb{N}$ and let $\mathcal{A}$ be an event in the stationary $\sigma$-field. Then $\{X_n\}$ conditioned on $\mathcal{A}$ is again a Markov chain, with transition probabilities*

$$\mathbb{P}(X_{n+1} = j | X_n = i, \mathcal{A}) = \frac{\mathbb{P}_j(\mathcal{A})}{\mathbb{P}_i(\mathcal{A})} \mathbb{P}(X_{n+1} = j | X_n = i),$$

*where $\mathbb{P}_i(\mathcal{A}) := \mathbb{P}(\mathcal{A} | X_0 = i)$ is supposed to be positive.*

PROOF.   Note that $\mathbb{P}(\mathcal{A} | X_{n+1} = j, X_n = i) = \mathbb{P}(\mathcal{A} | X_{n+1} = j) = \mathbb{P}_j(\mathcal{A})$ for any invariant event $\mathcal{A}$. Then

$$\mathbb{P}(X_{n+1} = j | X_n = i, \mathcal{A}) = \frac{\mathbb{P}(X_{n+1} = j, X_n = i, \mathcal{A})}{\mathbb{P}(X_n = i, \mathcal{A})}$$

$$= \frac{\mathbb{P}(\mathcal{A}|X_{n+1}=j, X_n=i)\mathbb{P}(X_{n+1}=j, X_n=i)}{\mathbb{P}(\mathcal{A}|X_n=i)\mathbb{P}(X_n=i)}$$

$$= \frac{\mathbb{P}_j(\mathcal{A})}{\mathbb{P}_i(\mathcal{A})}\mathbb{P}(X_{n+1}=j|X_n=i). \qquad \square$$

In our case, let $X_n$ be simple random walk with $X_0 = 1$, stopped at 0 and $h$. Then $\mathcal{A} := \{X_n = h \text{ eventually}\}$ is a stationary event and, as shown by a well-known martingale argument, $\mathbb{P}_i(\mathcal{A}) = i/h$. Therefore, the new transition probabilities in the first leg of the excursion are

$$p^{\text{there}}(i, i-1) = \frac{i-1}{2i}, \qquad p^{\text{there}}(i, i+1) = \frac{i+1}{2i} \qquad \text{for } i = 1, \ldots, h-1.$$

(Note the consistency property that these values do not depend on $h$.) Similarly, in the second leg, the transition probabilities are

$$p^{\text{back}}(j, j-1) = \frac{h+2-j}{2(h+1-j)},$$

$$p^{\text{back}}(j, j+1) = \frac{h-j}{2(h+1-j)} \qquad \text{for } j = h, h-1, \ldots, 1.$$

We will denote the measure given by these *there* and *back* Markov chains by $\mathbb{P}^{\text{there}}(\cdot)$ and $\mathbb{P}^{\text{back}}(\cdot)$, with corresponding expected values $\mathbb{E}^{\text{there}}(\cdot)$ and $\mathbb{E}^{\text{back}}(\cdot)$. If a chain is started at some value $i$, we denote the corresponding measure using a subscript $i$.

This consistency property of the transition probabilities was possibly first noted in [46], proving that these conditioned walks can be obtained by watching a three-dimensional Bessel process when it hits positive integer points. This also implies that the limit of our excursions $\mathcal{E}_h$, as $h \to \infty$, normalized to have height 1, is two Bessel first passage bridges put back-to-back, which is D. Williams' description of Itô's *Brownian excursion* conditioned to have height 1. For more details on the combinatorial aspects of Brownian excursion theory, see [47].

An observation similar to this consistency property is the following. Let $T_i$ denote the first hitting time of the value $i$ by the simple random walk. Then, conditioning a walk started at $v$ first on $\mathcal{A} = \{T_h < T_0\}$ and then further on $\mathcal{B} = \{T_i < T_j\}$, where $v \in [i, j] \subseteq [0, h]$, is the same as first conditioning on $\mathcal{B}$ and then, after hitting $i$, conditioning on $\mathcal{A}$. It follows that for any fixed integer $q \geq 1$, the distribution of the $q$th sub-up-excursion of type $i \nearrow j \searrow i$ in $\mathcal{E}_h$ (i.e., an interval of steps of $\mathcal{E}_h$ that is an up-excursion with minimum $i$ and maximum $j$), conditioned on its existence, is just the distribution of a standard up-excursion of height $j - i$. This observation shows that the large-$N$ limit of typical height-$h$ up-contours considered at the beginning of this section is indeed the up-contour given by two independent copies of $\mathcal{E}_h$.

For the computations in Section 5, we will often use the following lemma.

LEMMA 4.2.   *Let $T_i$ denote the first hitting time of the value $i$. Then*

$$\mathbb{P}_j^{\text{there}}(T_i < T_h) = \frac{(h-j)i}{(h-i)j} \quad and \quad \mathbb{P}_i^{\text{back}}(T_j < T_0) = \frac{(h+1-j)i}{(h+1-i)j}$$

*for $1 \le i \le j \le h$ in the first case and $0 \le i \le j \le h$ in the second.*

PROOF.   For a death-and-birth chain $\{X_n\}$ on $\mathbb{Z}^+$ with absorbing state 1 and transition probabilities $p_i = \mathbb{P}(X_{n+1} = i+1 | X_n = i)$ and $q_i = \mathbb{P}(X_{n+1} = i-1 | X_n = i)$, if we define

$$(4.1) \qquad\qquad \varphi(x) := \sum_{m=1}^{x-1} \prod_{i=2}^{m} \frac{q_i}{p_i}, \qquad \varphi(1) := 0,$$

then $\{\varphi(X_n)\}$ is a martingale as far as $X_n \in \{2, 3, \dots\}$ and

$$(4.2) \qquad\qquad \mathbb{P}_x(T_a < T_b) = \frac{\varphi(b) - \varphi(x)}{\varphi(b) - \varphi(a)}$$

for $1 \le a \le x \le b$; see [29], Section 5.3. For the chain $\mathbb{P}^{\text{there}}$, as computed after Lemma 4.1, we have $q_i/p_i = (i-1)/(i+1)$, so

$$(4.3) \qquad\qquad \varphi^{\text{there}}(x) = \sum_{m=1}^{x-1} \frac{2}{m(m+1)} = 2\left(1 - \frac{1}{x}\right)$$

and (4.2) then gives the first result. The second result follows from the relation $\mathbb{P}_j^{\text{back}}(T_i < T_0) = \mathbb{P}_{h+1-j}^{\text{there}}(T_{h+1-i} < T_{h+1})$ with the *there* chain understood as the first leg of $\mathcal{E}_{h+1}$.   □

The two independent copies of an excursion of height $h$ will be denoted by $\mathcal{E}_h$ and $\mathcal{E}_h^*$, and the cycle determined by them will be denoted by $C(\mathcal{E}_h, \mathcal{E}_h^*)$. Even though we understand the structure of the excursions $\mathcal{E}_h$ well, our Two Cautious Hikers algorithm suggests that computing the expected length $L(h) := \mathbb{E}|C(\mathcal{E}_h, \mathcal{E}_h^*)|$ will not be an easy task. Nevertheless, at least we can easily prove the following intuitively clear monotonicity result which will be used at the end of Section 5.

LEMMA 4.3.   *For any $h \in \mathbb{Z}^+$, $L(h) \le L(h+1)$.*

PROOF.   Consider an up-contour $C = C(\mathcal{E}_{h+1}, \mathcal{E}_{h+1}^*)$ and perform the following procedure to get an up-contour $C'$, distributed as $C(\mathcal{E}_h, \mathcal{E}_h^*)$.

Deleting all the maxima of the up-excursion $\mathcal{E}_{h+1}$ clearly gives an excursion of height $h$, denoted by $e$, which is distributed as $\mathcal{E}_h$. On the other hand, deleting the minima of an up-excursion can give several excursions of height at most $h$, at least one of which has exactly height $h$. Denote all of these excursions by $e_1, \ldots, e_k$, appearing in this order, and let us suppose that the height of $e_i$ is $h$, $1 \leq i \leq k$. Let the maximum of $e_1$ be $e_1(a) = m \leq \max e_i$ and the maximum of $e_k$ be $e_k(b) = n \leq \max e_i$. Let $a'$ be the first time when $e_i$ reaches $m$ and $b'$ be the last time when it is $n$. Now, build an excursion $f$ that starts like $e_1$, up to time $a$, then take the part of $e_i$ between $a'$ and $b'$, and finish with the part of $e_k$ from time $b$ until its end. It is easy to see that this $f$ has the distribution of $\mathcal{E}_h$.

The point of this construction is that taking the contour cycle $C'$ determined by $e$ and $f$ has the distribution of $C(\mathcal{E}_h, \mathcal{E}_h^*)$, while it is not difficult to see that each edge in $C'$ has a natural copy in $C$, thus $|C'| \leq |C|$. That is, we have coupled the random variables $C(\mathcal{E}_{h+1}, \mathcal{E}_{h+1}^*)$ and $C(\mathcal{E}_h, \mathcal{E}_h^*)$ so that $|C(\mathcal{E}_h, \mathcal{E}_h^*)| \leq |C(\mathcal{E}_{h+1}, \mathcal{E}_{h+1}^*)|$, which implies $L(h) \leq L(h+1)$.    $\square$

Instead of a "typical cycle," one could try to define the analogs of $\mathfrak{L}(n)$ and $L(h)$ using the cycle containing the origin. However, that would correspond to picking a cycle of diameter $n$ or of height $h$ from a large $N \times N$ box not uniformly, but weighted by the lengths of the cycles. Therefore, the expected length of the cycle of the origin, conditioned on having height $h$, is $\mathbb{E}(|C(\mathcal{E}_h, \mathcal{E}_h^*)|^2)/\mathbb{E}|C(\mathcal{E}_h, \mathcal{E}_h^*)|$, which is of the same order as $L(h)$ iff the second moment estimate $\mathbb{E}(|C(\mathcal{E}_h, \mathcal{E}_h^*)|^2) \leq K(\mathbb{E}|C(\mathcal{E}_h, \mathcal{E}_h^*)|)^2$ holds with some constant $K < \infty$. We conjecture that this is the case, but cannot prove it.

We will show in Section 5 that $L(h) \sim ch^{2\delta}$ with some constant $c$ and $\delta = (1 + \sqrt{17})/2$. How do we get the asymptotics for $\mathfrak{L}(n)$ from this? It is well known, and follows easily, for example, from our Lemma 4.2 or the Bessel process representation mentioned above, that the length of $\mathcal{E}_h$ divided by $h^2$ has a nontrivial limiting distribution as $h \to \infty$, with exponential lower and upper tails. Therefore,

$$(4.4) \qquad \mathbb{P}(h^2/K < \text{ diameter of } C(\mathcal{E}_h, \mathcal{E}_h^*) < Kh^2) > 1 - \exp(-cK).$$

Setting $K = C \log h$ and using the trivial polynomial bounds $n < \mathfrak{L}(n) < n^2$, from $L(h) = \Theta(h^{2\delta})$, we easily get

$$(4.5) \qquad \Omega(n^\delta / \log n) \leq \mathfrak{L}(n) \leq O(n^\delta \log n),$$

which gives the first part of Theorem 1.2.

On the other hand, with the identity

$$(4.6) \qquad \lim_{h \to \infty} \frac{\log P(h)}{\log h} = 2 \lim_{n \to \infty} \frac{\log \mathfrak{P}(n)}{\log n},$$

we have to be a bit more careful because now, we need the exponential tails not for $\mathcal{E}_h$ and $C(\mathcal{E}_h, \mathcal{E}_h^*)$, but for the cycle $C_o$ going through the origin.

However, the weighting by the length is only of polynomial order against the exponential decay in the tail:

$$\mathbb{P}(|C(\mathcal{E}_h, \mathcal{E}_h^*)| > Kh^4) \leq \mathbb{P}(\text{diameter of } C(\mathcal{E}_h, \mathcal{E}_h^*) > \sqrt{K}h^2) \leq \exp(-c\sqrt{K})$$

implies

$$(4.7) \qquad \mathbb{P}(h^2/K < \text{ diameter of } C_o < Kh^2) > 1 - Ch^4\exp(-cK).$$

This is slightly weaker than (4.4), but still good enough to imply (4.6), provided that at least one of those limits exists.

**5. The expected length of a typical contour.** Take a contour $C$ determined by a pair of compatible excursions of height $h$. Without loss of generality, we may assume that $C$ is an up-contour on level 0 and that the marginal excursions $X[a, a+2k]$ and $Y[b, b+2\ell]$ have endpoints $X_a = X_{a+2k} = 0$ and $Y_b = Y_{b+2\ell} = -h$. We may also assume that $h \geq 2$; otherwise, we would simply have $|C| = 4$. Besides $C$, there are many other contours on level 0 intersecting the interior of the smallest enclosing rectangle $R = [a+1, a+2k] \times [b+1, b+2\ell]$. By Lemma 3.1, all these contours are entire cycles, denoted by $C_1, \ldots, C_t$, given by compatible subexcursion pairs within $X[a, a+2k]$ and $Y[b, b+2\ell]$. Thus, if we denote by $T'$ the set of all horizontal edges in $R$ separating 0- and 1-level faces, and by $T''$ the set of such vertical edges, then the following recursion-type relation holds:

$$(5.1) \qquad |C| = |T'| + |T''| - \sum_{i=1}^{t} |C_i|.$$

This is the relation that we will turn into an actual recursion for $L(h)$, by taking expectations w.r.t. our measure on the independent pair of excursions $\mathcal{E}_h$ and $\mathcal{E}_h^*$.

A horizontal edge in $R$, between the faces $(n+1/2, m+1/2)$ and $(n+1/2, m+3/2)$, is in $T'$ if and only if there is an $i \in \{0, 1, \ldots, h\}$ such that $X_n = i$ and $\{Y_m, Y_{m+1}\} = \{-i, -i+1\}$. The condition for a vertical edge to be in $T''$ can be written analogously. Thus,

$$|T'| + |T''| = \sum_{i=0}^{h-1} V_{X[a,a+2k]}(i)U_{Y[b,b+2\ell]}(-i) + \sum_{i=0}^{h-1} U_{X[a,a+2k]}(i)V_{Y[b,b+2\ell]}(-i),$$

where $V_{X[a,a+2k]}(i)$ is the number of visits of $X[a, a+2k]$ to $i$, while $U_{X[a,a+2k]}(i)$ is the number of steps $n$ with $\{X_n, X_{n+1}\} = \{i, i+1\}$; the analogous numbers for $Y[b, b+2\ell]$ are $V_Y(i)$ and $U_Y(i)$. Taking expectation w.r.t. $(\mathcal{E}_h, \mathcal{E}_h^*)$, we get

$$(5.2) \qquad T(h) = T'(h) + T''(h) = 2\sum_{i=1}^{h} V_h(i)U_h(h-i),$$

where we have used the notation $T'(h) := \mathbb{E}|T'(\mathcal{E}_h)|$ and $V_h(i) := \mathbb{E}(V_{\mathcal{E}_h})(i)$, etc. and the fact that the variables $V_{\mathcal{E}_h}(i)$ and $U_{\mathcal{E}_h}(i)$ are independent of the analogous variables for $\mathcal{E}_h^*$.

Now, we will similarly translate the expected total length of the additional contour cycles $C_1, \ldots, C_t$ into terms of the pair $(\mathcal{E}_h, \mathcal{E}_h^*)$. The additional up-contours are given exactly by the sub-up-excursion pairs in which the subexcursion of $X[a, a+2k]$ is of the form $i \nearrow j \searrow i$ and the subexcursion of $Y[b, b+2\ell]$ is of the form $-j \nearrow -i \searrow -j$, where $1 \le i < j \le h - 1$. Similarly, the additional down-contours are given exactly by the sub-down-excursion pairs in which the $X$-sub-excursion is of the form $j \searrow i \nearrow j$ and the $Y$-sub-excursion is of the form $1 - i \searrow 1 - j \nearrow 1 - i$, where $1 \le i < j \le h$.

The observation about the subexcursions of $\mathcal{E}_h$ immediately preceding Lemma 4.2, the linearity of expected values and the multiplication of expected values for independent variables collectively imply that taking the expected value of the relation (5.1) w.r.t. our measure on $(\mathcal{E}_h, \mathcal{E}_h^*)$, we get the following recursion for $L(h)$:

$$(5.3) \quad \begin{aligned} L(h) = T(h) &- \sum_{1 \le i < j \le h-1} N_h(i,j) N_h(h-j, h-i) L(j-i) \\ &- \sum_{1 \le i < j \le h} M_h(i,j) M_h(h+1-j, h+1-i) L(j-i), \end{aligned}$$

where $N_h(i,j)$ is the expected number of $i \nearrow j \searrow i$ sub-up-excursions and $M_h(i,j)$ is the expected number of $j \searrow i \nearrow j$ sub-down-excursions of $\mathcal{E}_h$ for $0 \le i < j \le h$.

Therefore, we now need to compute the expected values $V_h(i)$, $U_h(i)$, $N_h(i,j)$ and $M_h(i,j)$ for our excursion measure $\mathcal{E}_h$. This will be done using Lemma 4.2.

LEMMA 5.1.  (i) *For the expected number of visits to $i$ by $\mathcal{E}_h$, we have*

$$V_h(i) = 2\frac{i(h-i)}{h} + 2\frac{i(h+1-i)}{h+1} \qquad \text{for } i = 1, \ldots, h-1,$$

*while $V_h(0) = 2$ and $V_h(h) = 2(1 - 1/(h+1))$.*

(ii) *For the expected number of $\{i, i+1\}$-crossings, we have*

$$U_h(i) = 2\frac{(h-i)(i+1)}{h} + 2\frac{(h-i+1)(i+1)}{h+1} - 2 \qquad \text{for } i = 0, \ldots, h-1.$$

PROOF.  Count the number of visits of the *there* leg to $i$, where $1 \le i \le h - 1$, using the following method. There will necessarily be a first visit to $i$. Then the next step is to $i - 1$ with probability $(i-1)/(2i)$ and to $i + 1$ with probability $(i+1)/(2i)$. In the first case, there will certainly be a next visit to $i$; in the second case, a second visit occurs with probability $\mathbb{P}^{\text{there}}_{i+1}(T_i < T_h) = \frac{(h-i-1)i}{(h-i)(i+1)}$, by Lemma 4.2. Adding up these possibilities, we find that a second visit to $i$ occurs with probability $1 - \frac{h}{2i(h-i)}$ and we reach $h$ without a second visit with probability $\frac{h}{2i(h-i)}$. After each visit to $i$, we have these same possibilities independently, hence the number of visits to $i$ is a geometric random variable with expected value $\frac{2i(h-i)}{h}$. Treating the *back* leg similarly, we have proven (i).

To prove (ii), first notice that

$$U_h(i) = \mathbb{E}^{\text{there}}(2|\{n : X_n = i, X_{n+1} = i+1\}| - 1)$$
$$+ \; \mathbb{E}^{\text{back}}(2|\{n : X_n = i+1, X_{n+1} = i\}| - 1).$$

Then the computation of these expectations is very similar to part (i). In the *there* leg, at the first visit to $i+1$, which is just after the first $(i, i+1)$-crossing, there will be a second $(i, i+1)$-crossing with probability $\mathbb{P}^{\text{there}}_{i+1}(T_i < T_h)$, and so on; thus, we have a geometric random variable with mean $1/(1 - \mathbb{P}^{\text{there}}_{i+1}(T_i < T_h)) = \frac{(h-i)(i+1)}{h}$. The *back* leg can be treated in the same way and so the proof is complete.  □

Turning to the computation of $N_h(i,j)$ and $M_h(i,j)$, note that $N_h(i,j) = N^{\text{there}}_h(i,j) + N^{\text{back}}_h(i,j)$ whenever $j < h$, and $N_h(i,h) = 1 + N^{\text{back}}_h(i,h)$, with the obvious notation. Similarly, $M_h(i,j) = M^{\text{there}}_h(i,j) + M^{\text{back}}_h(i,j)$ for $j < h$, and $M_h(i,h) = M^{\text{back}}_h(i,h)$.

LEMMA 5.2.   *For* $1 \le i < j \le h$, *we have*

$$N^{\text{there}}_h(i,j) = \frac{i}{h} \frac{h-i}{(j-i)(j+1-i)},$$

$$N^{\text{back}}_h(i,j) = \frac{i}{h+1} \frac{h+1-i}{(j-i)(j+1-i)},$$

$$M^{\text{there}}_h(i,j) = \frac{h-j}{h} \frac{j}{(j-i)(j+1-i)},$$

$$M^{\text{back}}_h(i,j) = \frac{h+1-j}{h+1} \frac{j}{(j-i)(j+1-i)}.$$

PROOF.   For $i < j$, let $R^{\text{there}}_h(j \searrow i)$ denote the expected number of segments in the trajectory of the *there* leg of $\mathcal{E}_h$ that go from $j$ to $i$, visiting only vertices in $\{i+1, i+2, \ldots, j-1\}$ meanwhile. Similarly, let $R^{\text{back}}_h(i \nearrow j)$ denote the expected number of segments in the trajectory of the *back* leg that go from $i$ to $j$, visiting only vertices in $\{i+1, i+2, \ldots, j-1\}$ meanwhile. Now, observe that

$$N^{\text{there}}_h(i,j) = R^{\text{there}}_h(j \searrow i) - R^{\text{there}}_h(j+1 \searrow i),$$

$$(5.4) \qquad N_h^{\text{back}}(i,j) = R_h^{\text{back}}(i \nearrow j) - R_h^{\text{back}}(i \nearrow j+1),$$

$$M_h^{\text{there}}(i,j) = R_h^{\text{there}}(j \searrow i) - R_h^{\text{there}}(j \searrow i-1),$$

$$M_h^{\text{back}}(i,j) = R_h^{\text{back}}(i \nearrow j) - R_h^{\text{back}}(i-1 \nearrow j).$$

Whenever the *there* leg is at $j$, either $\{T_h < T_i\}$ or $\{T_i < T_h\}$ will happen and in the latter case, there will be exactly one segment of the trajectory counted in $R_h^{\text{there}}(j \searrow i)$ before we return to $j$. Then we again have an independent try for $\{T_h < T_i\}$, and so on, until we actually reach $h$. Thus, the number of $j \searrow i$ crossings is a geometric random variable and

$$R_h^{\text{there}}(j \searrow i) = \frac{\mathbb{P}_j^{\text{there}}(T_i < T_h)}{1 - \mathbb{P}_j^{\text{there}}(T_i < T_h)} = \frac{(h-j)i}{h(j-i)},$$

where we have used Lemma 4.2 to get the second equality. Similarly,

$$R_h^{\text{back}}(i \nearrow j) = \frac{\mathbb{P}_i^{\text{back}}(T_j < T_0)}{1 - \mathbb{P}_i^{\text{back}}(T_j < T_0)} = \frac{(h+1-j)i}{(h+1)(j-i)}.$$

Plugging these results into the equations of (5.4), we get the identities of Lemma 5.2. □

Now, we plug the formulas of Lemma 5.1 into (5.2) and after simple algebraic manipulations, with the main tool being $\sum_{i=1}^{k} i^a = k^{a+1}/(a+1) + O(k^a)$, we get

$$(5.5) \qquad T(h) = \frac{16}{15} h^3 + O(h^2) \qquad \text{as } h \to \infty.$$

Next, collecting the coefficients of $L(m)$ for $m = j - i = 1, 2, \ldots, h-1$ in (5.3) into $Y(h, m)$, again with the same simple, but extensive, algebraic manipulations, we get

$$(5.6) \qquad L(h) = T(h) - \sum_{m=1}^{h-1} Y(h,m) L(m)$$

with

$$(5.7) \qquad Y(h,m) = \frac{4}{15}\left(\frac{h^3}{m^4} - \frac{m}{h^2}\right) - \frac{4}{3}\left(\frac{h}{m^2} - \frac{1}{m}\right)$$
$$+ \frac{O(h^2)}{m^4} + \frac{O(h)}{m^3} + \frac{O(1)}{m^2} + \frac{O(h^{-1})}{m} + O(h^{-2}),$$

where the implicit constants in $O(\cdot)$ are independent of both $h$ and $m$. Note that it is possible to get exact formulas in (5.5) and (5.7) using the lower order terms in the closed formula for $\sum_{i=1}^{k} i^a$. However, the use of symbolic computation software is strongly recommended here: the exact formulas fill a couple of pages, unfortunately.

It will be convenient to remove the $h^2$ factor from the denominator in (5.7), so we define $K(h) := h^2 L(h)$ and $Q(h) := h^2 T(h) = (16/15)h^5 + O(h^4)$ and

$$
\begin{aligned}
X(h, m) &:= \frac{h^2}{m^2} Y(h, m) \\
(5.8) \qquad &= \frac{4}{15}\left(\frac{h^5}{m^6} - \frac{1}{m}\right) - \frac{4}{3}\left(\frac{h^3}{m^4} - \frac{h^2}{m^3}\right) \\
&\quad + \frac{O(h^4)}{m^6} + \frac{O(h^3)}{m^5} + \frac{O(h^2)}{m^4} + \frac{O(h)}{m^3} + \frac{O(1)}{m^2},
\end{aligned}
$$

and thus rewrite (5.6) as

$$
(5.9) \qquad K(h) = Q(h) - \sum_{m=1}^{h-1} X(h, m) K(m).
$$

The exact explicit solution of such a recursion seems impossible, so how can it give the growth rate of $K(h)$? Before embarking on a search for the solution, we describe a naive approach which will show what the main ideas are.

Consider a simplified continuous analog of (5.9),

$$
(5.10) \qquad k(t) = \frac{16}{15}t^5 - \int_1^t \left[\frac{4}{15}\left(\frac{t^5}{s^6} - \frac{1}{s}\right) - \frac{4}{3}\left(\frac{t^3}{s^4} - \frac{t^2}{s^3}\right)\right] k(s)\, ds,
$$

which is a *linear Volterra integral equation of the second kind*. The general theory ([17]) says that since both $(16/15)t^5$ and the integral kernel are smooth in the domain $1 \le s \le t < \infty$, there exists a unique smooth solution $k(t)$. This $k(t)$ is actually the solution of the sixth order linear ordinary differential equation that we get by differentiating (5.10) six times:

$$
(5.11) \qquad 0 = k^{(6)}(t) + \frac{8}{t^4}k''(t) - \frac{32}{t^5}k'(t) + \frac{32}{t^6}k(t),
$$

with initial conditions $\{k(1), k'(1), \ldots, k^{(5)}(1)\}$ determined by the derivatives of $(16/15)t^5$ at $t = 1$. This ODE is a special case of *Euler's equation*

$$
t^n x^{(n)}(t) + a_1 t^{n-1} x^{(n-1)}(t) + \cdots + a_{n-1} t x'(t) + a_n x(t) = 0,
$$

where the $a_i$'s are constants. This equation has two singular points, 0 and $\infty$, both regular. See [16], Chapter 4, which will be our standard reference, for background. Now, with a change of variables $t = e^u$, we get a linear ODE with constant coefficients, which can be explicitly solved: a fundamental set of solutions is $\{t^\mu(\log t)^k\}$, where $\mu$ runs through the roots of the associated *indicial equation*,

$$(\mu)_n + (\mu)_{n-1}a_1 + \cdots + \mu a_{n-1} + a_n = 0,$$

with the notation $(\mu)_n := \mu(\mu - 1) \cdots (\mu - n + 1)$, while $k$ runs through the nonnegative integers less than the multiplicity of the root $\mu$. In the case of (5.11), the indicial equation is

$$(5.12) \quad \mu(\mu-1)(\mu-2)(\mu-3)(\mu-4)(\mu-5) + 8\mu(\mu-1) - 32\mu + 32 = 0,$$

with roots $\mu_{1,2} = 1$, $\mu_{3,4} = 4$ and $\mu_{5,6} = (5 \pm \sqrt{17})/2$, so we get that the general solution is

$$(5.13) \quad k(t) = c_1 t^{(\sqrt{17}+5)/2} + c_2 t^4 \log t + c_3 t^4 + c_4 t \log t + c_5 t + c_6 t^{(5-\sqrt{17})/2}.$$

With certain initial conditions, some coefficients $c_i$ in (5.13) might vanish, but we know that $ch^2 \leq L(h) \leq Ch^3$, so we expect $c't^4 \leq k(t) \leq C't^5$ and it seems likely that the growth rate will be $t^{(\sqrt{17}+5)/2}$ as $t \to \infty$.

If we use $Q(t)$ instead of $(16/15)t^5$, and the full $X(t,s)$ as the integral kernel, we arrive at a nonhomogeneous Euler-type ODE with nonconstant coefficients $a_i = a_i(t)$. These $a_i(t)$ and $Q^{(6)}(t)$ are still holomorphic at $t = \infty$, so we can use the *Frobenius method* (see [16], Section 4.8 and [42], Section 3.4), to find a fundamental set of solutions very similar to the constant coefficient case, with growth rates determined by the roots $\mu_1, \ldots, \mu_6$. Now, one possibility to make this argument rigorous would be to adapt the Frobenius method to the discrete equation (5.9), as $h \to \infty$, with differentiation replaced by the discrete difference operator $(\Delta K)(h) := K(h+1) - K(h)$. Although this adaptation seems possible, we will adopt a more traditional approach.

To begin the actual proof, consider the *generating function* $\kappa(z) := \sum_{m=1}^\infty K(m)z^m$. Because of the polynomial bounds on $K(m)$, the radius of convergence of $\kappa(z)$ around $z = 0$ is 1. The task is then to convert (5.9) into an ODE involving $\kappa(z)$, with a singularity at $z = 1$ and to read off the rate of growth of $K(m)$ from the size of the singularity of the solution $\kappa(z)$.

We multiply both sides of (5.9) by $z^h$ and sum the equation for all positive $h$'s. To write the resulting equation as an ODE, first notice that

$$(5.14) \qquad \sum_{h=2}^\infty \sum_{m=1}^{h-1} \frac{h^k}{m^\ell} K(m)z^h = \sum_{m=1}^\infty \frac{K(m)}{m^\ell} \sum_{n=m+1}^\infty n^k z^n.$$

It is not difficult to see that for $|z| < 1$,

$$\sum_{n=m+1}^\infty n^k z^n = \frac{z^{m+1}}{(1-z)^{k+1}} (k! f_k(z) + (k)_{k-1} f_{k-1}(z) m(1-z) + \cdots$$

$$(5.15) \qquad \qquad + k f_1(z) m^{k-1}(1-z)^{k-1} + m^k(1-z)^k),$$

where each $f_i(z)$ is an entire function with $f_i(1) = 1$. This formula implies, for example, that

$$(5.16) \qquad \rho(z) := \sum_{h=1}^{\infty} Q(h) z^h \stackrel{z=1}{\simeq} \frac{16/15}{(1-z)^6},$$

where the notation $f(z) \stackrel{z=z_0}{\simeq} g(z)$ means that $f(z) = g(z)(1 + \varepsilon(z))$, where $\varepsilon(z)$ is an entire function with $\varepsilon(z_0) = 0$.

Now, because of the $m^6$ denominator in $X(h, m)$, let us define $\phi(z) := \sum_{m=1}^{\infty} \frac{K(m)}{m^6} z^m$. Then $\sum_{m=1}^{\infty} \frac{K(m)}{m^5} z^m = \phi'(z) z$ and, in general,

$$(5.17) \quad \sum_{m=1}^{\infty} \frac{K(m)}{m^\ell} z^m = \phi^{(6-\ell)}(z) z^{6-\ell} + a_{\ell,5-\ell} \phi^{(5-\ell)}(z) z^{5-\ell} + \cdots + \phi'(z) z$$

for $0 \leq \ell \leq 5$, with $a_{\ell,j} \in \mathbb{Z}$ constants, $j = 2, \ldots, 5 - \ell$.

Plugging (5.15) into (5.14) and then using (5.17) and (5.16), we can turn (5.9) into the following ODE for $\phi(z)$, with a singularity at $z = 1$, and initial conditions at $z = 0$ easily computable from the first few values of $K(h)$:

$$\phi^{(6)}(z) z^6 + a_{0,5} \phi^{(5)}(z) z^5 + \cdots + \phi'(z) z$$

$$= \rho(z) + \phi^{(5)}(z) \frac{z^6}{1-z} \left( -\frac{4}{15} + \frac{4}{15} + \frac{4}{3} - \frac{4}{3} \right)$$

$$+ \phi^{(4)}(z) z^5 \left\{ \frac{f_1(z)}{(1-z)^2} \left( -\frac{4}{15}(5)_1 + \frac{4}{3}(3)_1 - \frac{4}{3}(2)_1 \right) + \frac{A_{4,0}}{1-z} \right\}$$

$$+ \phi^{(3)}(z) z^4 \left\{ \frac{f_2(z)}{(1-z)^3} \left( -\frac{4}{15}(5)_2 + \frac{4}{3}(3)_2 - \frac{4}{3}(2)_2 \right) \right.$$
$$\left. + \frac{f_1(z) A_{3,1}}{(1-z)^2} + \frac{A_{3,0}}{1-z} \right\}$$

$$+ \phi^{(2)}(z) z^3 \left\{ \frac{f_3(z)}{(1-z)^4} \left( -\frac{4}{15}(5)_3 + \frac{4}{3}(3)_3 \right) \right.$$
$$\left. + \frac{f_2(z) A_{2,2}}{(1-z)^3} + \frac{f_1(z) A_{2,1}}{(1-z)^2} + \frac{A_{2,0}}{1-z} \right\}$$

$$+ \phi'(z) z^2 \left\{ \frac{f_4(z)}{(1-z)^5} \left( -\frac{4}{15}(5)_4 \right) + \frac{f_3(z) A_{1,3}}{(1-z)^4} \right.$$
$$\left. + \frac{f_2(z) A_{1,2}}{(1-z)^3} + \frac{f_1(z) A_{1,1}}{(1-z)^2} + \frac{A_{1,0}}{1-z} \right\}$$

$$+ \phi(z) z \left\{ \frac{f_5(z)}{(1-z)^6} \left( -\frac{4}{15} 5! \right) + \frac{f_4(z) A_{0,4}}{(1-z)^5} + \cdots + \frac{A_{0,0}}{1-z} \right\},$$

where the $A_{i,j}$'s are constants, explicitly computable from the integers $a_{\ell,j}$ in (5.17) and the exact coefficients in the $O(\cdot)$ error terms in (5.8).

Now, changing variables $v = 1 - z$, that is, writing $\psi(v) := \phi(1-v)$ and $\sigma(v) := \rho(1-v)$, and multiplying both sides by $v^6$, we arrive at the equation

$$(5.18) \qquad v^6 \psi^{(6)}(v) \stackrel{v=0}{\simeq} \tfrac{16}{15} - 8\psi''(v)v^2 + 32\psi'(v)v - 32\psi(v).$$

The holomorphic correction factor on the right-hand side, as well as the initial conditions at $v = 1$, can be computed explicitly. This is again a nonhomogeneous, nonconstant coefficient version of (5.11), but we are now interested in the singularity of the solution at $v = 0$. [Notice that getting the same ODE is an accident, due only to the fact that the indicial equation (5.12) is invariant under the transformation $\mu \mapsto 5 - \mu$.] For the homogeneous version, the general Frobenius method ([42], Section 3.4) gives a fundamental set of solutions,

$$\psi_1(v) := f_{1,0}(v)v + f_{1,1}(v)v \log v + f_{1,2}(v)v \log^2 v,$$

$$\psi_2(v) := f_{2,0}(v)v + f_{2,1}(v)v \log v + f_{2,2}(v)v \log^2 v + f_{2,3}(v)v \log^3 v,$$

$(5.19)$

$$\psi_3(v) := f_{3,0}(v)v^4, \qquad\qquad \psi_4(v) := f_{4,0}(v)v^4 + f_{4,1}(v)v^4 \log v,$$

$$\psi_5(v) := f_{5,0}(v)v^{(5+\sqrt{17})/2}, \qquad \psi_6(v) := f_{6,0}(v)v^{(5-\sqrt{17})/2},$$

where the $f_{i,j}$'s are entire functions with $f_{i,0}(0) = 1$ and arbitrarily many coefficients in their Taylor expansions can be computed explicitly, the coefficients of which decay superexponentially. From (5.19), we get a solution of the nonhomogeneous equation by the standard method (see [16], Theorem 3.6.4):

$$(5.20) \qquad \psi(t) = \psi_h(t) + \sum_{i=1}^{6} \psi_i(t) \int_1^t \frac{W_i(s)}{W(s)} \sigma(s) s^6 \, ds,$$

where $s^6 \sigma(s) \stackrel{s=0}{\simeq} 16/15$ is the nonhomogeneity term, $W(t) := \det(\psi_i^{(j-1)}(t))_{i,j=1}^6$ is the Wronskian of the ODE (5.18), $W_i(t)$ is the determinant obtained from $W(t)$ by replacing the $i$th row $(\psi_i^{(j-1)}(t))_{j=1}^6$ by $(0,\ldots,0,1)$ and $\psi_h(t)$ is the solution of the homogeneous equation with the given initial conditions at $t = 1$. One can easily see that the integrand is an entire function, so the integral does not depend on the path of integration. It follows that the solution $\psi(t)$ of (5.18) has a leading term $\sim ct^{(5-\sqrt{17})/2}$ as $t \to 0$, unless the initial conditions accidentally kill this term. To show that this cancellation in (5.20) does not happen, it is enough to calculate the first few coefficients in the Taylor expansions of the $f_{i,j}(v)$'s in (5.19) because we can control the size of the error by the superexponential decay of these coefficients. We spare the reader these details.

From $\psi(v) \sim cv^{(5-\sqrt{17})/2}$, it follows that $\psi^{(4)}(v) \sim c'v^{(-3-\sqrt{17})/2}$ as $v \to 0$ and, hence, by (5.17),

$$\sum_{m=1}^{\infty} L(m)z^m = \sum_{m=1}^{\infty} \frac{K(m)}{m^2}z^m \sim c'(1-z)^{(-3-\sqrt{17})/2} \qquad \text{as } z \to 1.$$

The standard *Tauberian theorem for power series* ([30], Theorem 5 of Section XIII.5) gives that

$$\sum_{m=1}^{h} L(m) \sim c''h^{(3+\sqrt{17})/2}, \quad \text{moreover,} \quad L(h) \sim c'''h^{(1+\sqrt{17})/2} \qquad \text{as } h \to \infty,$$

where the second conclusion requires some additional hypothesis: for example, it is enough if $L(m)$ is monotone increasing. But we know this from Lemma 4.3, hence $L(h) \sim c'''h^{2\delta}$ with $2\delta = (1+\sqrt{17})/2$. Using this combined with the end of Section 4, we get (4.5) and the proof of the first part of Theorem 1.2 is complete.

## 6. The scaling relation $\gamma + \delta = 3/2$.

PROOF. Consider the $2N \times 2N$ box $B_N := \{(n,m) : \max\{|n|, |m|\} \le N\}$ around the origin. An edge $e = \{(n,m),(n+1,m)\}$, where $n+m$ is, say, even, is in a 0-level contour if $H(n+1/2, m+1/2) = 0$ and $H(n+1/2, m-1/2) = 1$. We will denote this event by $\mathcal{O}_e$. The classical De Moivre–Laplace theorem ([29], Section 2.1) tells us that $\mathbb{P}(X_n = k) = \Theta(n^{-1/2})$ for $k = O(\sqrt{n})$, which immediately implies that $\mathbb{P}(\mathcal{O}_e) = \Theta(N^{-1/2})$ for $(n,m) \in B_N$ and $\min\{|n|, |m|\} = \Omega(N)$, that is, for most edges inside the box $B_N$. [Here, the constant in the lower bound of $\Theta(N^{-1/2})$ depends, of course, on the constant in $\Omega(N)$.] From the linearity of expectation, it follows that the total expected length of 0-level contours inside $B_N$ is $\Theta(N^{3/2})$.

We call a closed contour a *medium cycle* if it goes through an edge inside $R := B_{2N/3} \setminus B_{N/3}$ and its height is between $h_*$ and $2h_*$, where $h_* = c_*\sqrt{N/\log N}$ with some small constant $c_* > 0$. Let us denote the total length of medium cycles contained entirely in $B_N$, a random variable, by $\lambda(N)$ and, if only those on some level $\ell \in \mathbb{Z}$ are continued, by $\lambda_\ell(N)$. Summing for all edges in $R$, by the linearity of expectation,

$$(6.1) \qquad \mathbb{E}\lambda(N) = \Theta(N^2(P(h_*) - P(2h_*))).$$

Note that if $\omega(N) := C\sqrt{\log N}$, with $C$ large enough, then the probability that there are contours in $B_N$ on any level $\ell$ with $|\ell| > \omega(N)\sqrt{N}$ is at most $o(N^{-1/2})$. Thus,

$$(6.2) \qquad \mathbb{E}\lambda(N) = \sum_{\ell \in \mathbb{Z}} \mathbb{E}\lambda_\ell(N) = \sum_{|\ell| < \omega(N)\sqrt{N}} \mathbb{E}\lambda_\ell(N) + o(N^{3/2}).$$

We now need the following simple lemma.

LEMMA 6.1.   *Let $\mathcal{S}$ be any finite set of faces and $F$ a face of $\mathbb{Z}^2$ such that $\tilde{H}(\mathcal{S} \cap F) \equiv 0$ is possible, where $\tilde{H}(n + \frac{1}{2}, m + \frac{1}{2}) := X_n + Y_m$. In other words, the graphical distance between any two faces in $\mathcal{S} \cup F$ is even. Then the conditional probability $\mathbb{P}(\tilde{H}(\mathcal{S}) \equiv 0 \mid \tilde{H}(F) = 2\ell)$, where $\ell \in \mathbb{Z}$, is maximized at $\ell = 0$. Furthermore, if $\mathrm{dist}(\mathcal{S}, F) > t$ and $\ell < \sqrt{t}/100$, then $\mathbb{P}(\tilde{H}(\mathcal{S}) \equiv 0 \mid \tilde{H}(F) = 2\ell) > \mathbb{P}(\tilde{H}(\mathcal{S}) \equiv 0 \mid \tilde{H}(F) = 0)/2$.*

PROOF.   We start by proving the first statement for the case when $|\mathcal{S}| = 1$. Suppose that $F = (1/2, 1/2)$ and take $n, m \in \mathbb{Z}$ such that $n + m$ is even. The stationarity and spatial reflection symmetry of SRW imply that

$$
\begin{aligned}
(6.3) \qquad &\mathbb{P}(X_n + Y_m = 0 | X_0 + Y_0 = 2\ell) \\
&= \mathbb{P}(X_n + Y_m = 2\ell | X_0 + Y_0 = 0) \\
&= \mathbb{P}(X_n + Y_m = 2\ell | X_0 = Y_0 = 0) \\
&= \sum_{j \in \mathbb{Z}} \mathbb{P}(X_n = j | X_0 = 0) \mathbb{P}(Y_m = j - 2\ell | Y_0 = 0).
\end{aligned}
$$

Now, it is well known that for any $n \in Z$, the sequences $\mathbb{P}(X_{2n} = 2\ell | X_0 = 0)$ and $\mathbb{P}(X_{2n+1} = 2\ell + 1 | X_0 = 0)$ are decreasing as $|\ell|$ increases and that $0 < a_1 \leq \cdots \leq a_k$ and $0 < b_1 \leq \cdots \leq b_k$ imply that $a_1 b_{\pi(1)} + \cdots + a_k b_{\pi(k)}$ is maximal when $\pi$ is the identical permutation. It follows that (6.3) is maximized at $\ell = 0$ and we are done. The case $|\mathcal{S}| > 1$ is similar—we just have to use the fact that for any finite set $\mathcal{N} \subset \mathbb{Z}$, the sequence $\mathbb{P}(X_{2n} = 2\ell$ for all $n \in \mathcal{N} | X_0 = 0)$ is decreasing in $|\ell|$.

For the second statement, first note that $\mathrm{dist}(\mathcal{S}, F) > t$ implies that there is a face $F^*$ sharing one coordinate with $F$, separating $F$ from $\mathcal{S}$ in the other coordinate, and such that $\mathrm{dist}(F^*, F) \geq t/2$. Recall, now, that two simple random walks started at distance $2\ell$ from each other with $\ell < \sqrt{t}/100$ can be coupled to coincide after at most $t/2$ steps with probability larger than $1/2$. This, the Markov property $\mathbb{P}(X_n = j | X_{t/2} = i, X_0 = 0) = \mathbb{P}(X_n = j | X_{t/2} = i)$ for any $n > t/2$ and the first statement of our lemma collectively imply the second statement.   $\square$

This lemma clearly implies that $\mathbb{E}\lambda_\ell(N) \leq \mathbb{E}\lambda_0(N)$ for all $\ell$ and that $\mathbb{E}\lambda_0(N) \leq 10\mathbb{E}\lambda_\ell(N)$ for $|\ell| < \sqrt{N}/100$. Hence, (6.2) gives

$$
\Omega(\sqrt{N}\mathbb{E}\lambda_0(N)) \leq \mathbb{E}\lambda(N) \leq \omega(N)\sqrt{N}\mathbb{E}\lambda_0(N) + o(N^{3/2}).
$$

This and (6.1) together give

$$
\begin{aligned}
(6.4) \qquad &\Omega\left( \frac{N^{3/2}}{\sqrt{\log N}} (P(h_*) - P(2h_*)) \right) - o(N) \\
&\leq \mathbb{E}\lambda_0(N) \leq O(N^{3/2}(P(h_*) - P(2h_*))).
\end{aligned}
$$

On the other hand, what is the number $\mu_0(N)$ of medium 0-level cycles in $B_N$? Consider just one box $D$ with side-length $N/3$ inside $R$. Combining the arguments of Section 2 and Lemma 6.1, we know that there is a medium 0-level cycle inside $D$ with a uniform positive probability, hence we have a lower bound $0 < m \leq \mathbb{E}\mu_0(N)$.

For an upper bound, consider the event $\mathcal{D} = \mathcal{D}(N, c^*, K) := \{$the diameter of each medium cycle is at most $Kh_*^2 \log h_*$, while the smaller sides of their enclosing rectangles are all at least $h_*^2/(K \log h_*)\}$. For any large $K$, if $c^*$ is small enough, then (4.7) implies that the probability of $\mathcal{D}$ is at least $1 - N^{-10}$, while $Kh_*^2 \log h_* \leq N/3$ and $h_*^2/(K \log h_*) \geq c'N/(K \log^2 N)$.

On the event $\mathcal{D}$, the maximum size of a set of medium cycles such that none of them contains another one in its interior is clearly bounded by $C_1 \log^4 N$ for some large constant $C_1$. On the other hand, the length of a nested sequence of 0-level medium cycles enclosing each other can be bounded as follows.

Fix a vertex $v$ and let $C_1, \ldots, C_t$ be the sequence of 0-level medium cycles whose enclosing rectangle contains $v$, with $C_i$ being in the interior of $C_{i+1}$, for all $i$. Note that if $C_i$ is an up-contour, then $C_{i+1}$ must be a down-contour and vice versa. By symmetry, we may assume that $C_1$ is an up-contour. Now, the down-excursions $e_2$, $e_2'$ that give $C_2$ must have endpoints strictly greater than the maximum of the up-excursions $e_1$, $e_1'$ giving $C_1$. This means that the SRW's $\{X_n\}$ and $\{Y_m\}$, in order to form $e_2$ and $e_2'$, respectively, must increase by at least the height of $C_1$ from the endpoints of $e_1$ and $e_1'$. Then, in order to form $e_3$ and $e_3'$, the up-excursions giving $C_3$, the SRW's must decrease by at least the height of $C_2$ from the endpoints of $e_2$ and $e_2'$, and so on. The time it takes for $\{X_n\}$ to form the sequence of up- and down-excursions $e_1, e_2, \ldots, e_t$ stochastically dominates the time it takes to form a similar sequence, but without the requirement that $e_i$ must be compatible with $e_i'$. Since $e_1$ has length at least $c'N/(K \log^2 N)$ and $e_t$ has length at most $N/3$, the expectation of $t$ is bounded by $C_2 \log^2 N$ for some large constant $C_2$.

Altogether, on $\mathcal{D}$, we have $\mu_0(N) \leq C_1 C_2 \log^6 N$. On $\mathcal{D}^c$, we still have $\mu_0(N) \leq N^2$, hence

$$(6.5) \qquad\qquad m \leq \mathbb{E}\mu_0(N) \leq C_3 \log^6 N.$$

The possible arrangement of all the 0-level cycles in these nested sequences is restricted by Lemma 3.1. Now, recall that the $j$th subexcursion of $\mathcal{E}_h$ with height $k < h$ has the distribution of $\mathcal{E}_k$. These two facts together imply that the cycles in the nested sequences can be chosen one by one so that each has the unconditional distribution of $\mathcal{E}_k$ for some $k$. This, our result (4.5) on $\mathfrak{L}(N)$ and (6.5) collectively imply

$$(6.6) \qquad\qquad \Omega(N^\delta/\log N) \leq \mathbb{E}\lambda_0(N) \leq O(N^\delta \log^7 N).$$

Now, comparing (6.4) with (6.6) yields

$$\Omega(N^{\delta-3/2}/\log N) \leq P(h_*) - P(2h_*) \leq O(N^{\delta-3/2}\log^{7.5} N) + o(N^{-1/2}).$$

Since $\delta > 1$, the term $o(N^{-1/2})$ is negligible and substituting back $N = \Theta(h_*^2 \log h_*)$ gives

$$\Omega(h_*^{2(\delta-3/2)}/\log^{1.5} h_*) \leq P(h_*) - P(2h_*) \leq O(h_*^{2(\delta-3/2)}\log^{7.5} h_*).$$

We can sum this over $h_* = 2^k h$ for $k = 0, 1, 2, \dots$ and any fixed $h$, arriving at

(6.7)       $$\Omega(h^{2(\delta-3/2)}/\log^{1.5} h) \leq P(h) \leq O(h^{2(\delta-3/2)}\log^{7.5} h).$$

That is, $P(h) \approx h^{-2\gamma}$ with $\gamma = 3/2 - \delta$. (4.6) then completes the proof of the second half of Theorem 1.2.   $\square$

## 7. Concluding remarks and open problems.

*The level sets of additive Brownian motion.*   In this subsection, we review the connections between our results and what is known about the additive Brownian motion.

Our Proposition 1.4 in the continuous setting is [25], Proposition 2.2. The main result of [24] concerning the additive Brownian motion is that given the zero set and the sign of a single excursion, the signs of all other excursions are determined. In corner percolation, if all 0-level curves are given, then deciding whether one cycle is an up- or down-contour is simply a matter of deciding which of the two possible chessboard colourings of the faces of $\mathbb{Z}^2$ to take and this obviously also determines the direction of all other cycles. Dalang and Mountford proved in [22] that there is a unique closed Jordan curve $J$ in the boundary $\partial B$ of the *Brownian bubble* $B$, the latter being defined as a connected component of the set $\{(s,t) \in \mathbb{R}^2 : B_t + B_s^* > 0\}$. It is natural to guess that the curves $C(\mathcal{E}_h, \mathcal{E}_h^*)$ have this $J$ as their scaling limit, and the Hausdorff dimension of $J$ is almost surely $(\sqrt{17}+1)/4$. Passing to the limit seems to be a nontrivial task, but, based on our results on the structure of the 0-level corner percolation set, we conjecture that the scaling limit of large cycles exists, is a simple loop and equals $J$. In [44], it was proven that $1 \leq \dim(J) \leq \dim(\partial B) < 3/2$ and the survey [19] announces the still unpublished result that $\dim(\partial B) = 3/2 - (5 - \sqrt{13 + 4\sqrt{5}})/4 = 1.421\dots$. The strict inequalities $\dim(J) < \dim(\partial B) < 3/2$ should be surprising only at first sight: on one hand, each macroscopic cycle $J$ has many tree-like sets attached to it, and the ones from the inside contribute significantly to $\partial B$; on the other hand, it was shown in [36] that almost all points of the level set are points of total disconnection, hence not in the boundary of any bubble.

Note that a Brownian excursion has a dense countable set of local extrema (see [45], Section 9) and two independent copies will not have any of them on the same height a.s. This means, on one hand, that there are no continuous analogs of the degree 4 points where the Two Cautious Hikers had to practice their caution and this lack of choice may serve as a simple intuitive explanation of the uniqueness of the Jordan curve inside $\partial B$. But, on the other hand, the denseness of the local extrema makes it unclear how to define a graph at all (which is why the definition of $J$ took up the entire paper [22]).

It could be interesting to study the set of outermost corner percolation cycles contained in a large finite box, and the scaling limit of this set. For conformally invariant models, such *loop ensembles* are the subject of active research; see [12, 13, 55]. One can also consider the tree structure of all the nested contour cycles. Could the scaling limit of this tree structure be described using Aldous' *continuum random tree* [1]?

*Biased coins.*  There are two natural versions of the model involving biased coins instead of fair ones. The first is when we define the $\pm$ sequences $\xi$ and $\eta$ using a $p$-biased coin for some $0 < p < 1$. Then the random walks $\{X_n\}$ and $\{Y_m\}$ will have every even step biased to one direction and every odd step biased to the other, so the main features of these walks are the same as if they were simple random walks. In particular, the proof of Theorem 1.1 in Section 2, relying on Brownian approximation, goes through verbatim. Furthermore, the formulas of Lemmas 4.2 and 5.1 will hold asymptotically for $h \to \infty$ and $i, j = \Theta(h)$, which means that the main coefficients in (5.9), hence the ODE (5.18), will be the same. Checking whether the cancellation of the leading term $\psi_6(t)$ in (5.20) happens for any value of $p$ would require some additional work which we have not done.

In the second version, the random walks $\{X_n\}$ and $\{Y_m\}$ are themselves biased to, say, the positive direction: $\mathbb{P}(X_{n+1} = X_n + 1) = \mathbb{P}(Y_{m+1} = Y_m + 1) = p > 1/2$. Here, one can easily prove that there are infinitely many infinite contour lines going from $(+\infty, -\infty)$ to $(-\infty, +\infty)$ with finite nested sequences of contour cycles between them. This model is no longer critical and the behavior of the sequences $\mathfrak{P}_p(n)$ and $\mathfrak{L}_p(n)$ is already a large deviation question. On the other hand, we can try to define the *near-critical exponent*

$$\mathbb{P}_p(\text{the contour of the origin is infinite}) = (p - 1/2)^{\beta + o(1)}.$$

It would be interesting to prove that this $\beta$ exists and to find its value.

*The cycle of the origin.*  Does the second moment estimate

$$\mathbb{E}(|C(\mathcal{E}_h, \mathcal{E}_h^*)|^2) \leq K(\mathbb{E}|C(\mathcal{E}_h, \mathcal{E}_h^*)|)^2$$

hold for some constant $K < \infty$? As explained in Section 4, this would imply that the cycle containing the origin, conditioned to have diameter $n$, has expected length $n^{\delta+o(1)}$. The analogous results are known for critical Bernoulli percolation; see [40].

*Crossing probability and corner percolation on a torus. Noise stability and dynamical corner percolation.* The criticality of Bernoulli(1/2) bond percolation on $\mathbb{Z}^2$ is intimately related to the fact that the *left-right crossing probability* in an $n \times n$ square is bounded away from 0 and 1; the exact value is actually 1/2 (with the appropriate choice of exactly what a left-right crossing means). In corner percolation, from the graph being 2-regular, it immediately follows that we cannot have both a left-right and an up-down crossing, while these events clearly have the same probability, hence this probability is at most 1/2. In a configuration without either type of crossing, the largest contours intersecting the sides of the square are all located on at most two levels, so our Lemma 3.1 can be applied to their possible arrangements. This does not leave many possibilities for noncrossing configurations, so we conjecture (also supported by computer simulations) that their probability tends to 0 as $n \to \infty$ and thus the left-right crossing probability tends to 1/2. A related result is that *corner percolation on the torus* $\mathbb{Z}_{2n} \times \mathbb{Z}_{2n}$ has a noncontractible cycle with probability tending to 1, or, equivalently, the infinite doubly-$2n$-periodic corner percolation configuration lifted to the universal cover $\mathbb{Z}^2$ will have an infinite path. This is because the only way to avoid this is to have $X_0 = X_{2n}$ and $Y_0 = Y_{2n}$ for the fundamental domain $[0, 2n] \times [0, 2n]$. However, we do not know the distribution of the homology classes of the noncontractible cycles.

The fact that the left-right crossing probability in critical Bernoulli percolation is 1/2 also gave rise to the study of *noise sensitivity* by Benjamini, Kalai and Schramm [9]: if we resample an arbitrary $\varepsilon > 0$ proportion of the percolation configuration in an $n \times n$ square, then having a left-right crossing in the new configuration will be asymptotically independent of having one previously. Moreover, Schramm and Steif [49] give a good estimate on the rate with which $\varepsilon = \varepsilon(n)$ can go to 0 so that we still have asymptotic independence, and this proves that critical percolation (at least on the triangular lattice) is also *dynamically sensitive*: if each variable in the configuration is updated independently according to a Poisson clock, then a.s. there will be exceptional times when there exists an infinite cluster. In corner percolation, the exceptional configurations with neither left-right nor up-down crossings are similar to those in which flipping a small number of signs changes the situation from having a left-right crossing to an up-down crossing, hence we conjecture that the event of a left-right crossing in corner percolation is *noise stable*: this event occuring in a configuration has an $n$-independent positive correlation with it occurring after each sign is resampled independently with probability $\varepsilon > 0$. Furthermore, the proof in Section 2 suggests that having no infinite components is a *dynamically stable* property, that is, there are no exceptional times a.s. In contrast, recurrence of simple random walk $\mathbb{Z}^2$ is dynamically sensitive [35].

*Is corner percolation almost supercritical?*   Benjamini asked if corner percolation is also critical from the point of view that adding to it an independent Bernoulli($\varepsilon$) bond percolation process results in a unique infinite cluster for any $\varepsilon > 0$. We conjecture that the answer is yes; however, it seems hard to construct an actual proof. The uniqueness part of this question for certain models with infinite components was shown in [8] and [10].

*Between the Gaussian free field and additive BM.*   An $\varepsilon$-almost-vertical *domino tiling* of the $n \times n$ square is a set of disjoint dimers that covers the $n \times n$ square, each dimer having at least one of its vertices inside this square and the proportion of the horizontal dimers being "close" to $\varepsilon > 0$. An $\varepsilon$-almost-horizontal domino tiling is defined analogously. Now, take a uniform random tiling from each family and take the union of them. One can prove that as $n \to \infty$, in the weak limit, we get a random 2-regular subgraph $G_\varepsilon$ of $\mathbb{Z}^2$. As with any double dimer model, this model has a natural height function. (Actually, our height function could also be defined by a version of this general procedure.) For the $\varepsilon$-almost-vertical tiling, if we scale the lattice by $1/n$ in both directions and do not scale the values of the height function, then the limit is a Gaussian free field $h_\varepsilon^\leftrightarrow$ with variance larger in the horizontal direction than in the vertical; see [37]. The same limit $h_\varepsilon$ of the height function for $G_\varepsilon$ is then the difference of a horizontally and a vertically stretched GFF: $h_\varepsilon = h_\varepsilon^\leftrightarrow - h_\varepsilon^\updownarrow$. As $\varepsilon \to 0$, can we rescale $h_\varepsilon$ so that it converges to additive BM? And what happens if one takes $\varepsilon$ tending to 0 simultaneously with $n$? Limits of such degenerate (single) dimer models appear in [11].

*The $k$-xor model on planar lattices.*   Benjamini suggested the following version of corner percolation on $\mathbb{Z}^2$: parametrize the edges $e$ of $\mathbb{Z}^2$ by their midpoints $(x_e, y_e) \in (\frac{1}{2}\mathbb{Z})^2$, with $x_e + y_e \in \mathbb{Z} + 1/2$. Now, take two sequences of i.i.d. Bernoulli($1/2$) variables, $\{\xi(k)\}$ and $\{\eta(k)\}$, parametrized by $k \in \frac{1}{2}\mathbb{Z}$, and let $\zeta(e) := \xi(x_e) + \eta(y_e) \bmod 2$. The resulting configuration of open ($\zeta = 1$) and closed ($\zeta = 0$) edges still has the property that two neighboring infinite lines either agree or complement each other, but the states of all the edges are now pairwise (moreover, 3-wise) independent. The connected components of this graph are no longer cycles, but one can show that the global behavior of the model is still governed by corner percolation, with the same critical exponents.

On the triangular lattice, to produce a critical model, we want site percolation with density 1/2. Benjamini's *trixor model* has three sequences of i.i.d. Bernoulli(1/2) variables, $\{\xi(k)\}$, $\{\eta(k)\}$ and $\{\zeta(k)\}$, parametrized by the three families of infinite parallel lines constituting the lattice. Then the state of a vertex in the triangular lattice is $\tau(v) := \xi(k) + \eta(\ell) + \zeta(j) \bmod 2$, given by the three lines through $v$. This is the same model as uniform measure on 0–1 configurations with the property that each vertex has an even number of 0-labeled (and hence an even number of 1-labeled) neighbours; see Figure 7. As observed by Angel and Schramm, trixor also has a natural height function, the sum of three independent simple random walks, with the components being the level sets. This means that the components and the contours separating them have similar long-range behavior, and our proof in Section 2 implies that there are only finite components a.s. It seems very likely that Section 6 can also be adapted to give $\gamma + \delta = 3/2$, provided that these exponents for trixor exist; namely, $\gamma$ would be the tail probability exponent for the cluster of the origin to reach out far, while $\delta$ would be the exponent for the expected length of a contour or (equivalently for trixor) the volume of a cluster. According to our simulations, these exponents are $\gamma \in (0.16, 0.2)$ and $\delta \in (1.3, 1.34)$. The exact combinatorial description of trixor contours using the three marginal simple random walks, in the spirit of our Proposition 1.4, seems to be harder than for corner percolation, and even if it is possible, the natural recursion will probably be in two variables, which then must be turned into a single variable recursion and then into an ODE. Finding $\gamma$ and $\delta$ for trixor would be a welcome development.

One can also take $k$ sequences of i.i.d. Bernoulli(1/2) variables, parametrized by $k$ families of suitable parallel lines, and then take the xor (i.e., mod 2 sum) of them to define a *$k$-xor* bond percolation on $\mathbb{Z}^2$ or site percolation on the triangular lattice. One might speculate that these models are all critical, with exponents $\gamma_k \to 5/48$ and $\delta_k \to 7/4$, which are the exponents for the $SLE_6$ process [43, 53]. But reality is even wilder: the simplest $k = 4$ choice on the triangular lattice, with the four families of lines having angles $0, 2\pi/3, 4\pi/3$ and $\pi/2$, already produced simulation results $\gamma \in (0.93, 1.05)$ and $\delta \in (1.74, 1.76)$. Further simulations by Braverman suggest that the hitting probabilities of the exploration path also coincide with those of $SLE_6$, supporting Conjecture 1.3 in the Introduction. Note that among the conformally invariant curves $SLE_\kappa$, the defining property of $\kappa = 6$ is locality [48, 54], a property that obviously holds for Bernoulli percolation, but presently we see no explanation as to how it arises for 4-xor.

As a variation on trixor, Angel suggested considering uniform measure on 0–1 configurations on the vertices of the triangular lattice, with the constraint that each vertex has an odd number of 1-labelled neighbours. This *odd-trixor* model is also given by flipping the states of a well-chosen deterministic quarter of the vertices in the trixor model; see Figure 8. Hence, it is again of linear entropy and simulations by Braverman suggest locality and $SLE_6$ hitting probabilities, again supporting Conjecture 1.3; see Figure 9. It would be very interesting to prove at least noise and dynamical sensitivity for 4-xor and odd-trixor. But, we no longer have a good height function interpretation, so we cannot even prove that all the components of 4-xor or odd-trixor are a.s. finite.

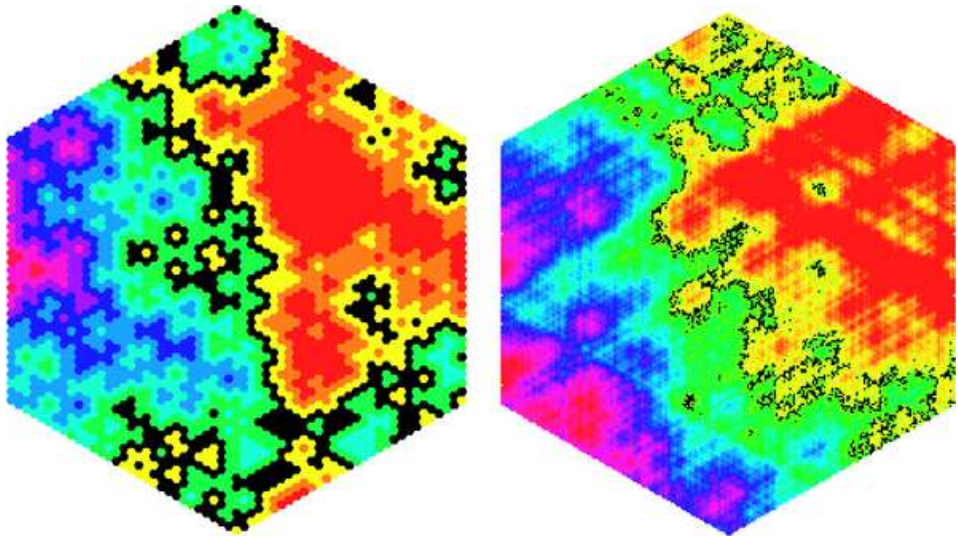

Fig. 7.   *Two trixor height function samples (side-lengths 30 and 200), with the 0-level set painted black.*

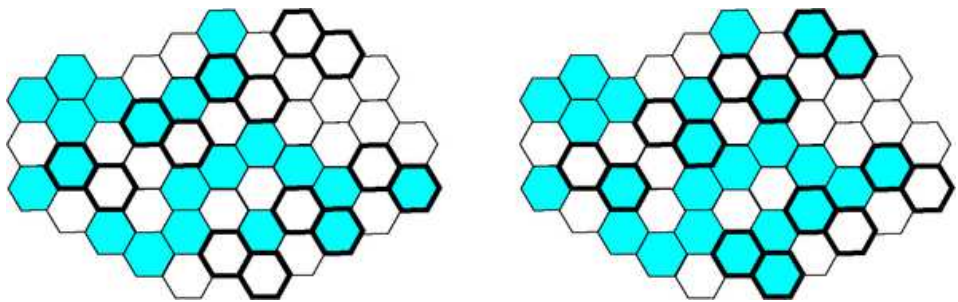

Fig. 8.   *From even- to odd-trixor, deterministically.*

As a final variation, Krishnapur suggested using other balanced symmetric Boolean functions of $k$ variables instead of $k$-xor (which is just parity). The simplest choice is majority, when $k$ is odd. We have not studied the resulting $k$-*majority* models, except for the observation that 3-majority has strong long-range dependence reflecting the triangular lattice.

As we learned from Diaconis [27], the natural 2-xor site percolation model on $\mathbb{Z}^2$ is used to disprove the conjecture in mathematical psychology that *human vision* cannot distinguish random patterns with the same first and second order statistics. Note that the methods of [27] break down if one wants to produce a random pattern that is visually different from Bernoulli percolation, but the same up to fourth order statistics (i.e., 4-wise independent). Trixor is 5-wise independent, but our models do not produce arbitrary large independence with a global behavior different from the one observed in Bernoulli percolation. Nevertheless, as the very recent work [7] shows, one might have up to roughly $\log n$-wise independence in an $n \times n$ square but a completely different global behavior.

There are many other ways to define loop models on the triangular or the hexagonal grid, mostly with quadratic entropy and conjectural conformal invariance; see [55, 56].

*Random walks on randomly oriented lattices.* Fix an arbitrary $\pm 1$ sequence $\{\xi(n)\}_{n \in \mathbb{Z}}$ and take a simple random walk $(X_j)_{j=0}^{\infty}$ on $\mathbb{Z}$ with $X_0 := 0$. Now, consider the two-dimensional random walk $(X_k, S_k)_{k=0}^{\infty}$ defined by $S_0 := 0$ and $S_k := \sum_{j=0}^{k-1} \xi(X_j)$. These processes were studied in [14, 15, 34].

Note that if $\xi$ is the i.i.d. random sequence $\mathbb{P}(\xi(n) = 1) = \mathbb{P}(\xi(n) = -1) = 1/2$, then $(X_k, S_k)_{k=0}^{\infty}$ is exactly the random walk we mentioned in the Introduction, interpolating between the almost simple random walk and the corner percolation component going through the origin.

Let us call a sequence $\xi$ *hospitable* if $(X_k, S_k) = 0$ infinitely often a.s., and *hostile* otherwise. Is there a characterization of hospitable sequences $\xi(n) \in \{\pm 1\}$? It was shown in [14] that a periodic sequence with the same number of $+1$'s and $-1$'s is hospitable, while the sequence $\xi(n) := \operatorname{sign}(n)$ and an i.i.d. fair coin random sequence are hostile a.s. The following is a rough explanation as to why.

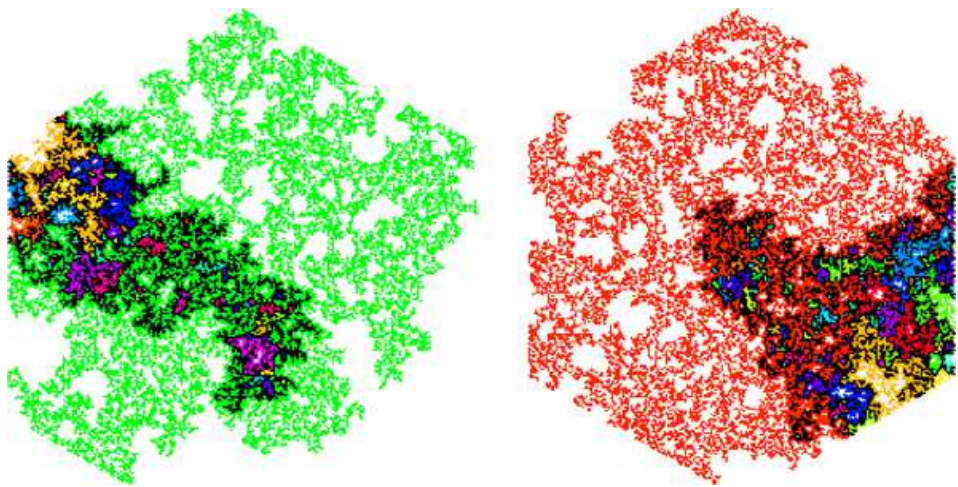

Fig. 9. *All the clusters neighboring the cluster of the origin (in black) in odd-trixor and in Bernoulli(1/2) site percolation on the triangular lattice.*

By the first $k$ steps, $\{X_j\}$ visits each element of an interval of the form $[-\Theta(\sqrt{k}), \Theta(\sqrt{k})]$ around $\Theta(\sqrt{k})$ times and almost all other elements of $\mathbb{Z}$ are visited much less. Therefore, if $|\sum_{n=-m}^{m} \xi(n)| = f(m)$ with $f(m)$ varying regularly in some sense, then we expect $S_k$ to grow like $\sqrt{k} f(\sqrt{k})$ and thus $\mathbb{P}((X_k, S_k) = (0,0))$ to decay like $1/(kf(\sqrt{k}))$. This is summable if $f(m) = \Omega(\log^{1+\varepsilon}(m))$ and then $\xi$ is hostile. Thus, a hospitable sequence must be quite uniform, without large fluctuations in the difference between the number of $+1$'s and $-1$'s. In particular, when $\xi$ is the i.i.d. fair coin sequence, then $f(m)$ is typically $\sqrt{m}$, hence $S_k$ should grow like $k^{3/4}$ (see [18] for precise results) and hence the return probabilities $\mathbb{P}((X_k, S_k) = (0,0))$ should decay like $k^{-5/4}$, which was actually proven in [14].

It is easy to see that hospitality is a shift-invariant property of $\xi$. It would be interesting to decide about invariance under finite permutations: if $\xi$ is hospitable and $\pi$ is a permutation of $\mathbb{Z}$ moving only finitely many elements, is $\xi \circ \pi = \{\xi(\pi(n))\}_{n \in \mathbb{Z}}$ then also hospitable? If $0 = \tau(0) < \tau(1) < \tau(2) < \ldots$ are the successive return times of $\{X_j\}$ to $0$ and $R_i := S_{\tau(i)} - S_{\tau(i-1)}$, then $\xi$ is hospitable iff the i.i.d. jump sequence $R_i$ gives a recurrent walk on $\mathbb{Z}$. It is not very difficult to show that a slightly stronger hypothesis, namely that $R_i$ satisfies the weak law of large numbers, is indeed invariant under finite permutations. Finally, hospitality is not invariant under wobbling bijections of $\mathbb{Z}$, that is, under permutations $\pi$ that move all elements of $\mathbb{Z}$ to at most a fixed bounded distance. For example, the alternating sequence $\xi(n) := (-1)^n$ is hospitable, but the sequence $\xi \circ \pi$ is hostile if $\pi(2n) = 2n - 1$ and $\pi(2n - 1) = 2n$ for $n \geq 1$, while $\pi(n) = n$ for $n \leq 0$.

**Acknowledgments.** I am indebted to Yuval Peres for helpful discussions and advice. In particular, having seen my original proof of Theorem 1.1, he suggested that a simpler proof must exist. I am grateful to Bálint Tóth for sharing the problem with me and to Omer Angel, Itai Benjamini, Mark Braverman and Oded Schramm for information on the trixor and odd-trixor models. Besides them, I would like to thank Cedric Boutillier, Alan Hammond, Ben Hough, Davar Khoshnevisan, Russell Lyons, Jim Pitman, Scott Sheffield and Ádám Timár for conversations and comments. Finally, I am indebted to a referee for important simplifications and corrections.

## REFERENCES

[1] Aldous, D. (1991). The continuum random tree. II. An overview. *Stochastic Analysis (Durham, 1990). London Math. Soc. Lecture Note Ser.* **167** 23–70. Cambridge Univ. Press. MR1166406

[2] BALISTER, P., BOLLOBÁS, B. and STACEY, A. (2000). Dependent percolation in two dimensions. *Probab. Theory Related Fields* **117** 495–513. MR1777130

[3] BAUER, M. and BERNARD, D. (2003). Conformal field theories of stochastic Loewner evolutions [CFTs of SLEs]. *Comm. Math. Phys.* **239** 493–521. MR2000927

[4] BAXTER, R. J. (1982). *Exactly Solved Models in Statistical Mechanics*. Academic Press, London. MR0690578

[5] BAXTER, R. J. (2002). Completeness of the Bethe Ansatz for the six and eight-vertex models. *J. Statist. Phys.* **108** 1–48. MR1909554

[6] BELAVIN, A., POLYAKOV, A. and ZAMOLODCHIKOV, A. (1984). Infinite conformal symmetry in two-dimensional quantum field theory. *Nuclear Phys. B* **241** 333–380. MR0757857

[7] BENJAMINI, I. GUREL-GUREVICH, O. and PELED, R. On k-wise independent events and percolation. In preparation.

[8] BENJAMINI, I., HÄGGSTRÖM, O. and SCHRAMM, O. (2000). On the effect of adding $\varepsilon$-Bernoulli percolation to everywhere percolating subgraphs of $\mathbb{Z}^d$. *J. Math. Phys.* **41** 1294–1297. MR1757959

[9] BENJAMINI, I., KALAI, G. and SCHRAMM, O. (1999). Noise sensitivity of Boolean functions and applications to percolation. *Inst. Hautes Études Sci. Publ. Math.* **90** 5–43. MR1813223

[10] BENJAMINI, I., LYONS, R., PERES, Y. and SCHRAMM, O. (2001). Uniform spanning forests. *Ann. Probab.* **29** 1–65. MR1825141

[11] BOUTILLIER, C. The bead model & limit behaviors of dimer models. Preprint. arXiv:math.PR/0607162.

[12] CAMIA, F., FONTES, L. R. G. and NEWMAN, C. M. (2006). Two-dimensional scaling limits via marked nonsimple loops. *Bull. Braz. Math. Soc.* **37** 537–559. MR2284886

[13] CAMIA, F. and NEWMAN, C. M. (2006). Two-dimensional critical percolation: The full scaling limit. *Comm. Math. Phys.* **268** 1–38. MR2249794

[14] CAMPANINO, M. and PÉTRITIS, D. (2003). Random walks on oriented lattices. *Markov Process. Related Fields* **9** 391–412. MR2028220

[15] CAMPANINO, M. and PÉTRITIS, D. (2004). On the physical relevance of random walks: An example of random walks on a randomly oriented lattice. In *Proceedings of Random Walks and Geometry, Vienna 2001* (V. Kaimanovich, ed.) 393–411. de Gruyter, Berlin. MR2087791

[16] CODDINGTON, E. A. and LEVINSON, N. (1955). *Theory of Ordinary Differential Equations*. McGraw-Hill, New York. MR0069338

[17] CORDUNEANU, C. (1991). *Integral Equations and Applications*. Cambridge Univ. Press. MR1109491

[18] CSÁKI, E., KÖNIG, W. and SHI, Z. (1999). An embedding for the Kesten–Spitzer random walk in random scenery. *Stochastic Process. Appl.* **82** 283–292. MR1700010

[19] DALANG, R. C. (2003). Level sets and excursions of the Brownian sheet. *CIME 2001 Summer School, Topics in Spatial Stochastic Processes. Lecture Notes in Math.* **1802** 167–208. Springer, Berlin. MR1975520

[20] DALANG, R. C. and MOUNTFORD, T. (1996). Nondifferentiability of curves on the Brownian sheet. *Ann. Probab.* **24** 182–195. MR1387631

[21] DALANG, R. C. and MOUNTFORD, T. (1997). Points of increase of the Brownian sheet. *Probab. Theory Related Fields* **108** 1–27. MR1452548

[22] DALANG, R. C. and MOUNTFORD, T. (2001). Jordan curves in the level sets of additive Brownian motion. *Trans. Amer. Math. Soc.* **353** 3531–3545. MR1837246

[23] DALANG, R. C. and MOUNTFORD, T. (2002). Eccentric behaviors of the Brownian sheet along lines. *Ann. Probab.* **30** 293–322. MR1894109

[24] DALANG, R. C. and MOUNTFORD, T. (2002). Nonindependence of excursions of the Brownian sheet and of additive Brownian motion. *Trans. Amer. Math. Soc.* **355** 967–985. MR1938741

[25] DALANG, R. C. and WALSH, J. B. (1993). The structure of a Brownian bubble. *Probab. Theory Related Fields* **96** 475–501. MR1234620

[26] DALANG, R. C. and WALSH, J. B. (1993). Geography of the level sets of the Brownian sheet. *Probab. Theory Related Fields* **96** 153–176. MR1227030

[27] DIACONIS, P. and FREEDMAN, D. (1981). On the statistics of vision: The Julesz conjecture. *J. Math. Psychology* **24** 112–138. MR0640207

[28] DOOB, J. L. (1959). Discrete potential theory and boundaries. *J. Math. Mech.* **8** 433–458; erratum 993. MR0107098

[29] DURRETT, R. (1996). *Probability: Theory and Examples*, 2nd ed. Duxbury Press, Belmont, CA. MR1609153

[30] FELLER, W. (1971). *An Introduction to Probability Theory and Its Applications.* II, 2nd ed. Wiley, New York. MR0270403

[31] FERRARI, P. L. and SPOHN, H. (2006). Domino tilings and the six-vertex model at its free fermion point. *J. Phys. A: Math. Gen.* **39** 10297–10306. MR2256593

[32] GÁCS, P. (2002). Clairvoyant scheduling of random walks. In *Proceedings of the Thirty-Fourth Annual ACM Symposium on Theory of Computing* 99–108. ACM, New York. Available at http://arxiv.org/abs/math.PR/0109152v7.

[33] GRIMMETT, G. (1999). *Percolation*, 2nd ed. Springer, Berlin. MR1707339

[34] GUILLOTIN-PLANTARD, N. and LE NY, A. Random walks on FKG-horizontally oriented lattices. Preprint. arXiv:math.PR/0303063.

[35] HOFFMAN, C. (2006). Recurrence of simple random walk on $\mathbb{Z}^2$ is dynamically sensitive. *ALEA Lat. Am. J. Probab. Math. Stat.* **1** 35–45. MR2235173

[36] KENDALL, W. (1980). Contours of Brownian processes with several-dimensional time. *Z. Wahrsch. Verw. Gebiete* **52** 269–276. MR0576887

[37] KENYON, R. (2001). Dominos and the Gaussian free field. *Ann. Probab.* **29** 1128–1137. MR1872739

[38] KENYON, R., OKOUNKOV, A. and SHEFFIELD, S. (2006). Dimers and amoebae. *Ann. Math.* **163** 1019–1056. MR2215138

[39] KESTEN, H. (1980). The critical probability of bond percolation on the square lattice equals $\frac{1}{2}$. *Comm. Math. Phys.* **74** 41–59. MR0575895

[40] KESTEN, H. (1986). The incipient infinite cluster in two-dimensional percolation. *Probab. Theory Related Fields* **73** 369–394. MR0859839

[41] KHOSHNEVISAN, D. and XIAO, Y. (2002). Level sets of additive Lévy processes. *Ann. Probab.* **30** 62–100. MR1894101

[42] KOHNO, M. (1999). *Global Analysis in Linear Differential Equations*. Kluwer, Dordrecht. MR1700776

[43] LAWLER, G. F., SCHRAMM, O. and WERNER, W. (2002). One-arm exponent for critical 2D percolation. *Electron. J. Probab.* **7** 1–13. MR1887622

[44] MOUNTFORD, T. S. (1993). Estimates of the Hausdorff dimension of the boundary of positive Brownian sheet components. *Séminaire de Probabilités XXVII. Lecture Notes in Math.* **1557** 233–255. Springer, Berlin. MR1308568

[45] PERES, Y. (2001). *An Invitation to Sample Paths of Brownian Motion*. Lecture notes at UC Berkeley. Available at http://www.stat.berkeley.edu/~peres/.

[46] PITMAN, J. (1975). One-dimensional Brownian motion and the three-dimensional Bessel process. *Adv. in Appl. Probab.* **7** 511–526. MR0375485

[47] PITMAN, J. (2006). *Combinatorial Stochastic Processes. Lecture Notes in Math.* **1875**. MR2245368

[48] SCHRAMM, O. (2000). Scaling limits of loop-erased random walks and uniform spanning trees. *Israel J. Math.* **118** 221–288. MR1776084

[49] SCHRAMM, O. and STEIF, J. (2008). Quantitative noise sensitivity and exceptional times for percolation. *Ann. Math.* Available at http://arxiv.org/abs/math.PR/0504586.

[50] SHEFFIELD, S. (2005). Random surfaces. *Astérisque* **304**. MR2251117

[51] SHEFFIELD, S. (2007). Gaussian free fields for mathematicians. *Probab. Theory Related Fields* **139** 521–541. MR2322706

[52] SMIRNOV, S. (2001). Critical percolation in the plane: Conformal invariance, Cardy's formula, scaling limits. *C. R. Acad. Sci. Paris Sér. I Math.* **333** 239–244. Expanded version available at http://www.math.kth.se/~stas/papers/index.html. MR1851632

[53] SMIRNOV, S. and WERNER, W. (2001). Critical exponents for two-dimensional percolation. *Math. Res. Lett.* **8** 729–744. MR1879816

[54] WERNER, W. (2004). Random planar curves and Schramm–Loewner evolutions. *Lectures Notes on Probability Theory and Statistics. Lecture Notes in Math.* **1840** 107–195. Springer, Berlin. MR2079672

[55] WERNER, W. (2005). Some recent aspects of random conformally invariant systems. Preprint. arXiv:math.PR/0511268.

[56] WILSON, D. (2004). On the Red–Green–Blue model. *Phys. Rev. E* **69** 037105.

[57] WINKLER, P. (2000). Dependent percolation and colliding random walks. *Random Structures Algorithms* **16** 58–84. MR1728353

MICROSOFT RESEARCH
ONE MICROSOFT WAY
REDMOND, WASHINGTON 98052-6399
USA
E-MAIL: gabor@microsoft.com
URL: http://research.microsoft.com/~gabor/